\documentclass[11pt,a4paper]{article}

\usepackage[hyperfootnotes=false,pdfpage layout=SinglePage]{hyperref}
\usepackage{mathptmx}
\usepackage[scaled=.92]{helvet}
\usepackage{amssymb}
\usepackage{amsmath}
\usepackage{amsthm}
\usepackage{latexsym}
\usepackage{amsfonts}
\usepackage{mathrsfs}
\usepackage{setspace}
\usepackage{graphicx}
\usepackage{stmaryrd}
\usepackage{pgf,tikz}

\newtheorem{thm}{Theorem}[section]
\newtheorem{cor}[thm]{Corollary}
\newtheorem{lem}[thm]{Lemma}
\newtheorem{prop}[thm]{Proposition}
\theoremstyle{definition}
\newtheorem{defn}[thm]{Definition}
\newtheorem{rem}[thm]{Remark}

\numberwithin{equation}{section}
\renewcommand{\emptyset}{\varnothing}

\newcommand{\R}{\ensuremath{\mathbb R}}    
\newcommand{\C}{\ensuremath{\mathbb C}}    
\newcommand{\N}{\ensuremath{\mathbb N}}    

\newcommand{\gperp}{{[\perp]}}

\newcommand{\product}{[\cdot\,,\cdot]}
\newcommand{\hproduct}{(\cdot\,,\cdot)}

\newcommand{\calB}{\mathcal B}         \newcommand{\frakB}{\mathfrak B}
         
\newcommand{\calD}{\mathcal D}

\newcommand{\calG}{\mathcal G}         
\newcommand{\calH}{\mathcal H}         
\newcommand{\calI}{\mathcal I}         
\newcommand{\calJ}{\mathcal J}         \newcommand{\frakJ}{\mathfrak J}
         
\newcommand{\calL}{\mathcal L}         
\newcommand{\calM}{\mathcal M}

\newcommand{\calR}{\mathcal R}

\newcommand{\calU}{\mathcal U}

\newcommand{\la}{\lambda}
\newcommand{\veps}{\varepsilon}
\newcommand{\vphi}{\varphi}

\newcommand{\mat}[4]
{
   \begin{pmatrix}
      #1 & #2\\
      #3 & #4
   \end{pmatrix}
}
\newcommand{\vek}[2]
{
   \begin{pmatrix}
      #1\\
      #2
   \end{pmatrix}
}

\renewcommand{\Im}{\operatorname{Im}}
\renewcommand{\Re}{\operatorname{Re}}

\renewcommand{\ker}{\operatorname{ker}}
\newcommand{\ran}{\operatorname{ran}}
\newcommand{\dom}{\operatorname{dom}}

\newcommand{\sap}{\sigma_{{ap}}}

\renewcommand{\sp}{\sigma_{+}}
\newcommand{\sm}{\sigma_{-}}

\newcommand{\Lra}{\Longrightarrow}

\newcommand{\Llra}{\Longleftrightarrow}

\newcommand{\ol}{\overline}
\newcommand{\ds}{\dotplus}
\newcommand{\wt}{\widetilde}

\newcommand{\sgn}{\operatorname{sgn}}
\newcommand{\dist}{\operatorname{dist}}
\newcommand{\AC}{\operatorname{AC}}

\newcommand{\fraka}{\mathfrak{a}}
\newcommand{\frakt}{\mathfrak{t}}
\DeclareMathOperator*{\liim}{l.i.m.}

\begin{document}
\vspace*{-.3cm}
\begin{center}
\begin{spacing}{1.7}
{\LARGE\bf Indefinite Sturm-Liouville operators with periodic coefficients}
\end{spacing}

\vspace{.5cm}
\begin{spacing}{1.7}
{\Large Friedrich Philipp}
\end{spacing}
\end{center}

\vspace{.1cm}\hrule\vspace*{.4cm}
\noindent{\bf Abstract}

\vspace*{.3cm}\noindent
We investigate the spectral properties of the maximal operator $A$ associated with a differential expression $\frac 1 w\left(-\frac d {dx}\left(p\frac d {dx}\right) + q\right)$ with real-valued periodic coefficients $w$, $p$ and $q$ where $w$ changes sign. It turns out that the non-real spectrum of $A$ is bounded, symmetric with respect to the real axis and consists of a finite number of analytic curves. The real spectrum is band-shaped and neither bounded from above nor from below. We characterize the finite spectral singularities of $A$ and prove that there is only a finite number of them. Finally, we provide a condition on the coefficients which ensures that $\infty$ is not a spectral singularity of $A$.
\vspace*{.4cm}\hrule
%
%

\vspace*{.5cm}
\section*{Introduction}
In this paper we consider the maximal operator $A$ corresponding to the Sturm-Liou\-ville differential expression
\begin{equation}\label{e:a}
\fraka(f) := \frac 1 w\Big(\!-(pf')' + qf\Big)
\end{equation}
on $\R$ with real-valued coefficients $w$, $p$ and $q$ which are periodic with the period $a > 0$ such that $w$, $q$ and $p^{-1}$ are integrable over $(0,a)$. We assume $p > 0$ and $w\neq 0$ almost everywhere.

It is well-known (see, e.g., \cite{w3}) that in the definite case (i.e.\ $w > 0$ (a.e.)) the operator $A$ is self-adjoint in the weighted $L^2$-Hilbert space $L^2_w(\R)$. Moreover, it is bounded from below and its spectrum has a band structure, i.e.\ it consists of compact intervals which may intersect in their endpoints only.

If the weight function $w$ changes its sign, the differential expression $\fraka$ in \eqref{e:a} is called {\it indefinite}. In this case the operator $A$ is no longer self-adjoint in $L^2_{|w|}(\R)$. But if $J$ denotes the operator of multiplication with $\sgn(w(\cdot))$, then $JA$ is the maximal operator associated with the definite Sturm-Liouville expression
\begin{equation}\label{e:t}
\frakt(f) := \frac 1{|w|}\Big(\!-(pf')' + qf\Big)
\end{equation}
and is therefore self-adjoint in the Hilbert space $L^2_{|w|}(\R)$. The operator $A$ is a so-called $J$-self-adjoint operator in $L^2_{|w|}(\R)$. Equivalently, $A$ is self-adjoint with respect to the {\it indefinite} inner product $\product := (J\cdot,\cdot)$, where $\hproduct$ denotes the scalar product in $L^2_{|w|}(\R)$.

If the lower bound of the spectrum of the self-adjoint operator $JA$ is positive, the spectral properties of the operator $A$ are well understood (see \cite{dl,ko,mz}). Such problems are called {\it left-definite}. In this case the spectrum of $A$ is real and has a band structure as in the definite case. But it is neither bounded from above nor from below if $w$ is indefinite. Also in the non-periodic case the $J$-self-adjoint operators corresponding to left-definite problems have real spectra and have been intensively studied in the literature. Here, we only mention \cite{be,bb,bbw,bv,kwz1,kwz2} and the monograph \cite{z} for further references. If $A_0$ is a $J$-self-adjoint realization of $\fraka$ and only $JA_0\ge 0$ is assumed, the spectrum of the operator $A_0$ remains real provided its resolvent set is non-empty. In this case it is of particular interest whether the operator $A_0$ is similar to a self-adjoint operator. With regard to this problem we refer to \cite{kk,kkm,km} and also to \cite{ko} in the periodic case. If the negative spectrum of $JA_0$ only consists of a finite number of negative eigenvalues (counting multiplicities) and $\rho(A_0)\neq\emptyset$, the spectrum of the operator $A_0$ is real with the possible exception of a finite number of non-real eigenvalues, cf.\ \cite{bt,bm,cl,l}. The situation becomes much more difficult, in general, when the negative spectrum of the operator $JA_0$ has accumulation points in $(-\infty,0]$. For example, accumulation of the non-real spectrum of $A_0$ to the real line may occur. In \cite{b,bp} and \cite{kt} such problems have been tackled with the local spectral theory of self-adjoint operators in Krein spaces.

In the present paper we discuss the periodic case and allow the lower bound of the spectrum of the operator $JA$ to be negative. Our methods are based on the interplay between standard tools in the analysis of periodic ODEs (such as Gelfand transform and Floquet-discriminant) and elements of the local spectral theory of $J$-self-adjoint operators. As a first result we prove that the Floquet-discriminant of $\fraka$ in \eqref{e:a} is not a constant. This implies that the spectrum of $A$ consists of the closures of analytic curves. Moreover, it is symmetric with respect to $\R$ since $A$ is $J$-self-adjoint. But in contrast to the left-definite case the non-real part of the spectrum might be non-empty. However, the main result of this paper (see Theorem \ref{t:dt}) shows that the non-real spectrum of $A$ is bounded. Moreover, Theorem \ref{t:dt} states that there is a finite number of real points which separate the real axis into open intervals of positive or negative type. This means that the spectrum of $A$ within these intervals is separable and the spectral subspaces corresponding to closed subintervals are Hilbert spaces with respect to either $\product = (J\cdot,\cdot)$ or $-\product$. Roughly speaking, the operator $A$ acts locally like a self-adjoint operator in a Hilbert space. The statements of Theorem \ref{t:dt} particularly imply that there is only a finite number of spectral singularities of $A$. These are characterized in Theorem \ref{t:reg}. In the last result (see Theorem \ref{t:infty}) a condition on the coefficients $w$ and $p$ of $\fraka$ is presented which ensures that the point $\infty$ is not a spectral singularity of $A$.

The paper is organized as follows. In section \ref{s:jsa} the necessary definitions and statements concerning $J$-self-adjoint operators are provided. In section \ref{s:mult_floq} we define the maximal operator $A$ and the multiplication operator $\wt A$ with the family of $J$-self-adjoint operators $\{A(t) : t\in [-\pi,\pi]\}$, where $A(t)$ is the differential operator in $L^2_{|w|}(0,a)$ associated with $\fraka$ on $(0,a)$ subject to the boundary conditions
$$
f(a) = e^{it}f(0)\quad\text{and}\quad (pf')(a) = e^{it}(pf')(0).
$$
The multiplication operator $\wt A$ acts in the Hilbert space $L^2([-\pi,\pi],L^2_{|w|}(0,a))$ and is unitarily equivalent to $A$. 
Since the spectrum of each operator $A(t)$ is non-empty (Proposition \ref{p:cl}), the Floquet-discriminant $D$ of $\fraka$ is a non-constant entire function. This implies that the spectrum of $A$, which coincides with $\{\la : D(\la)\in [-2,2]\}$, consists of closures of analytic curves, cf.\ Theorem \ref{t:A_first}. In section \ref{s:closer} we prove our main result Theorem \ref{t:dt}. The proof is based on perturbation-theoretic arguments. In section \ref{s:sf} it is proved that $A$ possesses a spectral function with a finite number of singularities. These are among the points in the spectrum of $A$ in which the derivative of $D$ vanishes. In Theorem \ref{t:reg} it is shown in particular that such a point $\la_0\in\sigma(A(t_0))$ is not a spectral singularity of $A$ if and only if the root subspace of $A(t_0)$ corresponding to $\la_0$ coincides with $\ker(A(t_0) - \la_0)$. The behaviour of the spectral function at $\infty$ is investigated in section \ref{s:infty}. We prove that $\infty$ is not a spectral singularity of $A$ if the weight function $w$ has only a finite number of turning points in $(0,a)$ and if $w$ and $p$ satisfy some regularity conditions in neighborhoods of these turning points.

\section{Preliminaries on \texorpdfstring{$J$}{J}-selfadjoint operators}\label{s:jsa}
In this paper $\calB(X,Y)$ denotes the set of all bounded and everywhere defined linear operators from a Banach space $X$ to a Banach space $Y$. As usual, we write $\calB(X)$ instead of $\calB(X,X)$. Spectrum and resolvent set of a closed linear operator $T$ are denoted by $\sigma(T)$ and $\rho(T)$, respectively. Point spectrum, residual spectrum and continuous spectrum of $T$ are denoted by $\sigma_p(T)$, $\sigma_r(T)$ and $\sigma_c(T)$, respectively. We set $\R^+ := (0,\infty)$ and $\R^- := (-\infty,0)$. Moreover, $\C^+$ ($\C^-$) denotes the upper (lower) complex halfplane.

Throughout this section let $(\calH,\hproduct)$ be a Hilbert space and let $J\in\calB(\calH)$ be boundedly invertible such that
$$
J = J^{-1} = J^*,
$$
where $J^*$ denotes the adjoint of $J$ with respect to the scalar product $\hproduct$. Such an operator will be called a {\it fundamental symmetry} in $\calH$. The fundamental symmetry $J$ induces a second inner product
$$
[f,g] := (Jf,g),\quad f,g\in\calH,
$$
on $\calH$ which is indefinite unless $J = \pm I$. The inner product space $(\calH,\product)$ is called a {\it Krein space}. For an intensive study of Krein spaces and operators therein we refer to the monographs \cite{ai} and \cite{b}.

A linear operator $B$ in $\calH$ is called {\it $J$-self-adjoint} if the operator $JB$ is self-adjoint. Equivalently, $B^* = JBJ$, which implies that the spectrum of $B$ is symmetric with respect to the real axis, i.e.
$$
\sigma(B) = \{\ol\la : \la\in\sigma(B)\}.
$$
Moreover, for $\la\in\C$ the following holds:
\begin{equation}\label{e:ran_ker}
\ker(B - \la) = J\ker(J(B - \la)J) = J\ker(B^* - \la) = J\ran(B - \ol\la)^\perp.
\end{equation}
Note that an operator is $J$-self-adjoint if and only if it is self-adjoint with respect to the inner product $\product$.

For the rest of this section let $B$ be a $J$-self-adjoint operator in $\calH$. Recall that the {\it approximate point spectrum} $\sap(B)$ of $B$ is defined as the set of all $\la\in\C$ for which there exists a sequence $(f_n)\subset\dom B$ with $\|f_n\| = 1$ and $(B - \la)f_n\to 0$ as $n\to\infty$. A point $\la\in\C$ is not an element of $\sap(B)$ if and only if $\ran(B - \la)$ is closed and $\ker(B - \la) = \{0\}$. Therefore, in view of \eqref{e:ran_ker} we have
\begin{equation}\label{e:R_sap}
\sigma(B)\cap\R\subset\sap(B).
\end{equation}
It should be mentioned that in general the spectral properties of $J$-self-adjoint operators differ considerably from those of self-adjoint operators. There exist simple examples of $J$-self-adjoint operators whose spectrum covers the entire complex plane or is empty. Therefore, the existing literature mainly focusses on special classes of $J$-self-adjoint operators such as definitizable or fundamentally reducible operators. Another approach is based on the local spectral analysis of $J$-self-adjoint operators.

\begin{defn}\label{d:spp}
Let $B$ be a $J$-self-adjoint operator in $\calH$. A point $\la\in\sap(B)$ is called a {\em spectral point of positive {\rm (}negative{\rm )} type} of $B$ if for every sequence $(f_n)\subset\dom B$ with $\|f_n\| = 1$ and $\|(B - \la)f_n\|\to 0$ as $n\to\infty$ we have
$$
\liminf_{n\to\infty}\,[f_n,f_n] > 0\quad\Big(\limsup_{n\to\infty}\,[f_n,f_n] < 0,\;\text{respectively}\Big).
$$
The set of all spectral points of positive {\rm (}negative{\rm )} type of $B$ will be denoted by $\sp(B)$ {\rm (}$\sm(B)$, respectively{\rm )}. A set $\Delta\subset\R$ is said to be of positive {\rm (}negative{\rm )} type with respect to $B$ if
$$
\Delta\cap\sigma(B)\subset\sp(B)\quad\Big(\Delta\cap\sigma(B)\subset\sm(B),\,\text{respectively}\Big).
$$
The set $\Delta$ is said to be of {\em definite} type with respect to $B$ if it is either of positive or negative type with respect to $B$.
\end{defn}

In the following we collect a few properties of the spectra of definite type. Proofs of the statements can be found in the fundamental paper \cite{lmm}. First of all we note that the spectra of positive and negative type of a $J$-self-adjoint operator $B$ are real. An isolated eigenvalue $\la$ of $B$ with finite (algebraic) multiplicity is of positive (negative) type if and only if the inner product $\product$ is positive definite (negative definite, respectively) on $\ker(B - \la)$. If $\calJ\subset\R$ is an interval which is of positive (negative) type with respect to $B$, then there exists an open domain $\calU$ in $\C$ such that $\calJ\subset\calU$ and $\calU\cap\sigma(B)\subset\sp(B)$ ($\calU\cap\sigma(B)\subset\sm(B)$, respectively). In particular, $\sp(B)$ and $\sm(B)$ are open in $\sigma(B)$. Moreover, the operator $B$ has a local spectral function $E$ on $\calJ$.

\begin{defn}\label{d:lsf}
Let $\Xi\subset\C$ be Borel-measurable and let $T$ be a closed linear operator in a Banach space $X$. By $\frakB_0(\Xi)$ we denote the system of Borel-measurable subsets of $\Xi$ whose closure is contained in $\Xi$. A mapping $E : \frakB_0(\Xi)\to\calB(X)$ is called a {\em local spectral function of $T$ on $\Xi$} if it has the following properties ($\Delta\in\frakB_0(\Xi)$):
\begin{enumerate}
\item[{\rm (S1)}] $E(\Delta)$ is a projection in the double-commutant of the resolvent of $T$ {\rm (}which maps to a {\rm (}closed{\rm )} subspace of $\dom T$ if $\Delta$ is bounded{\rm )}.
\item[{\rm (S2)}] $E$ is strongly $\sigma$-additive, i.e., if $\Delta_1,\Delta_2,\ldots\in\frakB_0(\Xi)$ are mutually disjoint and $\bigcup_{k=1}^\infty\,\Delta_k\in\frakB_0(\Xi)$, then
$$
\lim_{n\to\infty}\,\left\|E\left(\bigcup_{k=1}^\infty\,\Delta_k\right)f - \sum_{k=1}^n\,E(\Delta_k)f\right\| = 0
$$
holds for every $f\in X$.
\item[{\rm (S3)}] $E(\Delta_1\cap\Delta_2) = E(\Delta_1)E(\Delta_2)$ for all $\Delta_1,\Delta_2\in\frakB_0(\Xi)$.
\item[{\rm (S4)}] $\sigma(T|E(\Delta)\calH)\,\subset\,\ol{\sigma(T)\cap\Delta}$.
\item[{\rm (S5)}] $\sigma(T|(I-E(\Delta))\calH)\,\subset\,\ol{\sigma(T)\setminus\Delta}$.
\end{enumerate}
\end{defn}

The local spectral function $E$ of the $J$-self-adjoint operator $B$ on the interval $\calJ$ (which is of positive (negative) type with respect to $B$) has the additional property that the spectral subspace $E(\Delta)\calH$ is a Hilbert space with respect to the inner product $\product$ ($-\product$, respectively) for each $\Delta\in\frakB_0(\calJ)$. Hence, since $B$ is self-adjoint with respect to $\product$, the restriction $B|E(\Delta)\calH$ is self-adjoint in the Hilbert space $(E(\Delta)\calH,\product)$.

In the following, the signature of an inner product $\langle\cdot,\cdot\rangle$ on a finite-dimensional subspace $\calM$ of $\calH$ (which is contained in the domain of $\langle\cdot,\cdot\rangle$) will be denoted by
$$
\big\{\kappa_+(\langle\cdot,\cdot\rangle,\calM),\kappa_-(\langle\cdot,\cdot\rangle,\calM),\kappa_0(\langle\cdot,\cdot\rangle,\calM)\big\}.
$$
For $\la\in\C$ denote the root subspace of $B$ corresponding to $\la$ by
$$
\calL_\la(B) := \bigcup_{k=1}^\infty\ker\big((B - \la)^k\big).
$$
It is well-known (see, e.g., \cite[Proposition 3.2]{ds}) that for $\la\neq\ol\mu$ we have
$$
\calL_\la(B)\;\gperp\;\calL_\mu(B),
$$
i.e.\ $[f,g] = 0$ for all $f\in\calL_\la(B)$ and all $g\in\calL_\mu(B)$. In particular, $[f,g] = 0$ holds for all $f,g\in\calL_\la(B)$ if $\la\notin\R$, or equivalently, $\kappa_+(\product,\calL_\la(B)) = \kappa_-(\product,\calL_\la(B)) = 0$. Moreover, it is well-known that $\dim\calL_\la(B) = \dim\calL_{\ol\la}(B)$ holds for isolated eigenvalues $\la$ of $B$, see, e.g., \cite[Proposition I.3.2]{l}. The implication
\begin{equation}\label{e:sp_ker}
\la\in\sp(B)\cup\sm(B)\quad\Lra\quad\calL_\la(B) = \ker(B - \la)
\end{equation}
follows with the use of the local spectral function but also by elementary means: assume $(B - \la)f_1 = f_0$ and $(B - \la)f_0 = 0$, $f_0\neq 0$. Then $[f_0,f_0] = [(B - \la)f_1,f_0] = [f_1,(B - \la)f_0] = 0$ which contradicts $\la\in\sp(B)\cup\sm(B)$.

The next lemma describes the relation between the signatures of the inner products $\product$ and $[B\cdot,\cdot]$ on the subspace
$$
\calM_\la(B) := \calL_\la(B) + \calL_{\ol\la}(B),\quad\la\in\C.
$$

\begin{lem}\label{l:signatures}
Let $\la\in\C$ be an isolated eigenvalue of $B$ with finite multiplicity. Then the following holds. 
\begin{enumerate}
\item[{\rm (i)}]   If $\la\in\R$, then we have $\kappa_0(\product,\calL_\la(B)) = 0$. If, in addition, $\la\neq 0$, then also $\kappa_0([B\cdot,\cdot],\calL_\la(B)) = 0$.
\item[{\rm (ii)}]  If $\la\in\R^+$, then
$$
\kappa_\pm([B\cdot,\cdot],\calL_\la(B)) = \kappa_\pm(\product,\calL_\la(B)).
$$
\item[{\rm (iii)}] If $\la\in\R^-$, then
$$
\kappa_\pm([B\cdot,\cdot],\calL_\la(B)) = \kappa_\mp(\product,\calL_\la(B)).
$$
\item[{\rm (iv)}]  If $\la\in\C\setminus\R$, then also $\ol\la$ is a pole of the resolvent of $B$ with finite algebraic multiplicity and the Jordan structures of $B$ at $\la$ and $\ol\la$ coincide. Moreover,
$$
\kappa_-(\product,\calM_\la(B)) = \kappa_+(\product,\calM_\la(B)) = \dim\calL_\la(B).
$$
The same holds with respect to the inner product $[B\cdot,\cdot]$:
$$
\kappa_-([B\cdot,\cdot],\calM_\la(B)) = \kappa_+([B\cdot,\cdot],\calM_\la(B)) = \dim\calL_\la(B).
$$
In particular, $\dim\calL_\la(B) = \dim\calL_{\ol\la}(B)$ and
$$
\kappa_0(\product,\calM_\la(B)) = \kappa_0([B\cdot,\cdot],\calM_\la(B)) = 0.
$$
\end{enumerate}
\end{lem}
\begin{proof}
If $P(\la)$ denotes the Riesz-Dunford spectral projection of $B$ corresponding to $\la$, then $[P(\la)f,g] = [f,P(\ol\la)g]$ for $f,g\in\calH$, cf.\ \cite[Proposition 3.2]{l}. In particular, $P(\la)$ is $J$-self-adjoint if $\la$ is real. In this case we have
$$
\calH = \calL_\la(B)\,[\ds]\,(I - P(\la))\calH,
$$
where $[\ds]$ denotes the $\product$-orthogonal direct sum. Hence, if $f\in\calL_\la(B)$ with $[f,g] = 0$ for all $g\in\calL_\la(B)$, then $[f,g] = 0$ for all $g\in\calH$ and $f=0$ follows. This proves $\kappa_0(\product,\calL_\la(B)) = 0$ for real $\la$ and also $\kappa_0([B\cdot,\cdot],\calL_\la(B)) = 0$ if $\la\in\R\setminus\{0\}$. For the proof of (ii) we may assume that $\dim\calH < \infty$ and that $\sigma(B) = \{\la\}$. Using the Riesz-Dunford calculus, we define a square root $B^{1/2}$ of $B$. The operator $B^{1/2}$ is boundedly invertible and $J$-self-adjoint. Therefore, (ii) follows from $[Bf,f] = [B^{1/2}f,B^{1/2}f]$, $f\in\calH$.
The statement (iii) is proved similarly with the difference that $iB^{1/2}$ is $J$-self-adjoint. Statement (iv) is a consequence of \cite[Proposition 3.2]{l}.
\end{proof}

\begin{rem}
If the origin belongs to the spectrum of $B$, then there is in general no relation between the signatures of $\product$ and $[B\cdot,\cdot]$ on $\calL_0(B)$.
\end{rem}

The $J$-self-adjoint operator $B$ is said to have $\kappa$, $\kappa\in\N_0$, negative squares, if the inner product $[B\cdot,\cdot]$ has $\kappa$ negative squares. Equivalently (as $[B\cdot,\cdot] = (JB\cdot,\cdot)$), the spectrum of the self-adjoint operator $JB$ in $(-\infty,0)$ consists of exactly $\kappa$ eigenvalues (counting multiplicities). It should be mentioned that a $J$-self-adjoint operator with $\kappa$ negative squares and non-empty resolvent set is definitizable in the sense of \cite{l}. If $B$ has $\kappa = 0$ negative squares (i.e.\ $JB$ is non-negative or $[Bf,f]\ge 0$ for all $f\in\dom B$), then $B$ is called {\it $J$-nonnegative}.

\begin{prop}\label{p:kns}
Assume that the resolvent set of the $J$-self-adjoint operator $B$ is non-empty and that the resolvent of $B$ is compact. If $B$ has $\kappa$ negative squares, then
\begin{equation}\label{e:equ_kappa}
\sum_{\la\in\C^+\cup\R}\!\!\kappa_-\big([B\cdot,\cdot],\calM_\la(B)\big) \,=\,\kappa.
\end{equation}
In particular {\rm (}cf.\ Lemma {\rm\ref{l:signatures}}{\rm )}, the number of non-real eigenvalues of $B$ {\rm (}counting multiplicities{\rm )} does not exceed $\kappa$.
\end{prop}
\begin{proof}
Let $\la_1,\ldots,\la_n$ be the distinct non-zero eigenvalues of $B$ in $\C^+\cup\R$ with the property $\kappa_-([B\cdot,\cdot],\calM_{\la_j}(B)) > 0$, $j=1,\ldots,n$, and set
$$
\calM := \calL_0(B)\,[\ds]\,\calM_{\la_1}(B)\,[\ds]\,\ldots\,[\ds]\,\calM_{\la_n}(B).
$$
Then $\kappa_-([B\cdot,\cdot],\calM) = \sum_{\la\in\C^+\cup\R}\kappa_-([B\cdot,\cdot],\calM_\la(B))$. Since $\calM$ is $B$-invariant, the same holds for $\calM^\gperp = J\calM^\perp$, and we have $\calH = \calM[\ds]\calM^\gperp$ (see \cite[Theorem I.5.2]{l}). Hence, it remains to show that $[Bf,f]\ge 0$ for all $f\in\dom B\cap\calM^\gperp$. To see this, note that the spectrum of $B|\calM^\gperp$ is real and that $\R^+$ ($\R^-$) is of positive (negative) type with respect to $B|\calM^\gperp$, cf.\ Lemma \ref{l:signatures}. The same holds for the compact $J$-self-adjoint operator $C := (B|\calM^\gperp)^{-1}$. Hence, due to \cite[Corollary II.5.3]{l} we have $[Cf,f]\ge 0$ for all $f\in\calM^\gperp$, which proves the assertion.
\end{proof}

For $J$-self-adjoint operators $B$ with $\kappa$ negative squares and compact resolvent we set
$$
\sigma_{\rm ex}(B) := \big\{\la\in\C : \kappa_-\big([B\cdot,\cdot],\calM_\la(B)\big) > 0\big\}.
$$
The points in $\sigma_{\rm ex}(B)$ will be called the {\it exceptional eigenvalues} of $B$. It follows from Proposition \ref{p:kns} that $B$ has at most $\kappa$ exceptional eigenvalues in $\C^+\cup\R$ and hence a total of at most $2\kappa$ exceptional eigenvalues. The assertions of the next lemma follow directly from Lemma \ref{l:signatures} and \eqref{e:sp_ker}.

\begin{lem}\label{l:sigma_ex_basic}
Let $B$ be a $J$-self-adjoint operator with compact resolvent and $\kappa$ negative squares. Then a point $\la\in\sigma(B)\setminus\{0\}$ is contained in $\sigma_{\rm ex}(B)$ if and only if one of the following holds:
\begin{enumerate}
\item[{\rm (a)}] $\la\notin\R$.
\item[{\rm (b)}] $\la > 0$ and $\la\notin\sp(B)$.
\item[{\rm (c)}] $\la < 0$ and $\la\notin\sm(B)$.
\end{enumerate}
If $0\in\sp(B)\cup\sm(B)$, then $0\notin\sigma_{\rm ex}(B)$.
\end{lem}
%

\section{Multiplication operators and Floquet theory}\label{s:mult_floq}
The object of investigation in this paper will be the maximal operator associated with a Sturm-Liouville expression of the form
\begin{equation}\label{e:fraka}
\fraka(f) := \frac 1 w\Big(\!-(pf')' + qf\Big)
\end{equation}
on $\R$ with real-valued coefficients $w$, $p$ and $q$ which are periodic with the same period $a > 0$. We assume that the functions $w$, $q$ and $p^{-1}$ are integrable over $(0,a)$, that $w(x)\ne 0$ and $p(x) > 0$ for a.e.\ $x\in (0,a)$. If neither $w > 0$ a.e.\ nor $w < 0$ a.e.\ on $(0,a)$, we say that the weight function $w$ and the differential expression $\fraka$ are {\it indefinite}. Since also $|w|$ is $a$-periodic, it is a well-known fact (see, e.g., \cite[Lemma 12.1]{w3}) that the (definite) differential expression
\begin{equation}\label{e:frakt}
\frakt(f) := \frac 1 {|w|}\Big(\!-(pf')' + qf\Big)
\end{equation}
is limit point at $\pm\infty$. Hence, the maximal operator $T$ associated with $\frakt$ is self-adjoint in the weighted $L^2$-space $L^2_{|w|}(\R)$ which consists of all (equivalence classes of) measurable functions $f : \R\to\C$ such that $f^2w$ is integrable over $\R$. The scalar product on $L^2_{|w|}(\R)$ is given by
$$
(f,g) := \int_\R\,f(x)\ol{g(x)}|w(x)|\,dx,\quad f,g\in L^2_{|w|}(\R),
$$
and the maximal operator $T$ associated with $\frakt$ is defined by $Tf := \frakt(f)$ for $f\in\dom T$, where
$$
\dom T := \left\{f\in L^2_{|w|}(\R) : f,pf'\in\AC_{\rm loc}(\R),\,\frakt(f)\in L^2_{|w|}(\R)\right\}.
$$
Hereby, the set of all (locally) absolutely continuous complex-valued functions, defined on a bounded or unbounded interval $\Delta$, is denoted by $\AC(\Delta)$ ($\AC_{\rm loc}(\Delta)$, respectively). The maximal operator $A$ associated with $\fraka$ is defined analogously:
$$
Af := \fraka(f),\quad f\in\dom A := \dom T.
$$
Obviously, we have $JA = T$, where $J$ is the operator of multiplication with $\sgn(w(\cdot))$:
$$
(Jf)(x) := \sgn(w(x))f(x),\quad f\in L^2_{|w|}(\R),\;x\in\R.
$$
Since $J$ is a fundamental symmetry in $L^2_{|w|}(\R)$, the operator $A$ is $J$-self-adjoint. Equivalently, $A$ is self-adjoint with respect to the (in general indefinite) Krein space inner product
$$
[f,g] := (Jf,g) = \int_\R\,f(x)\ol{g(x)}w(x)\,dx,\quad f,g\in L^2_{|w|}(\R).
$$
The spectral properties of the operator $A$ are closely connected to those of a family $\{A(t) : t\in [-\pi,\pi]\}$ of differential operators associated with $\fraka$ in the Hilbert space $L^2_{|w|}(0,a)$. By $\hproduct_a$ we denote the scalar product in this Hilbert space, i.e.
$$
(f,g)_a := \int_0^a\,f(x)\ol{g(x)}|w(x)|\,dx,\quad f,g\in L^2_{|w|}(0,a).
$$
The operators $A(z)$, $z\in\C$, are defined by $A(z)f = \fraka(f)$ for $f\in\dom A(z)$, where
\begin{align*}
\dom A(z) := \{f\in L^2_{|w|}(0,a) : \,
&f,pf'\in\AC([0,a]),\;\fraka(f)\in L^2_{|w|}(0,a),\\
&f(a) = e^{iz}f(0),\;(pf')(a) = e^{iz}(pf')(0)\}.
\end{align*}
Note that $A(z+2\pi) = A(z)$ for all $z\in\C$. The operator of multiplication $J_a$ with the restriction of the function $\sgn(w(\cdot))$ to $[0,a]$ is a fundamental symmetry in the Hilbert space $L^2_{|w|}(0,a)$ and the operators $T(t) := J_aA(t)$ are self-adjoint in $L^2_{|w|}(0,a)$ for $t\in\R$. Hence, each of the operators $A(t)$, $t\in\R$, is $J_a$-self-adjoint and thus self-adjoint with respect to the inner product $\product_a$, where
$$
[f,g]_a := (J_af,g)_a = \int_0^a\,f(x)\ol{g(x)}w(x)\,dx,\quad f,g\in L^2_{|w|}(0,a).
$$
It is well-known that each operator $T(t)$, $t\in\R$, has compact resolvent and is bounded from below. By $\kappa(t)$ we denote the number of negative eigenvalues of $T(t)$ (counting multiplicities). Hence, the operator $A(t)$ has $\kappa(t)$ negative squares. Since each operator $T(t)$ has as many negative eigenvalues as either $T(0)$ or $T(\pi)$, see \cite[Theorem 12.7]{w3}, we have
$$
\kappa(t)\in\{\kappa^*-1,\kappa^*\},
$$
for all $t\in\R$, where $\kappa^*$ denotes the maximum of the number of negative eigenvalues of the operators $T(0)$ and $T(\pi)$.

\begin{prop}\label{p:cl}
The resolvent set of each operator $A(t)$, $t\in\R$, is non-empty and the resolvent of $A(t)$ is compact. In particular, we have
\begin{equation}\label{e:sigsum}
\sum_{\la\in\C^+\cup\R}\!\!\kappa_-\big([A(t)\cdot,\cdot]_a,\calM_\la(A(t))\big) \,=\,\kappa(t).
\end{equation}
Moreover, if the weight function $w$ is indefinite, then the real spectrum of $A(t)$ is neither bounded from below nor from above.
\end{prop}
\begin{proof}
The first two assertions are due to \cite[Corollary 1.4 and Proposition 2.2]{cl}, and \eqref{e:sigsum} follows from Proposition \ref{p:kns}. Assume, e.g., that $\sigma(A(t))\cap\R$ is bounded from below and choose $c > 0$ such that $-c < \min(\sigma(A(t))\cap\R)$ and $\max(\sigma_{\rm ex}(A(t))\cap\R) < c$. Then $[c,\infty)$ is of positive type with respect to $A(t)$ and $(-\infty,-c]\subset\rho(A(t))$. Therefore, see \cite[page 39]{l}, the point $\infty$ is not a singularity of the spectral function $E$ of $A(t)$, and it follows that $(E(\R\setminus [-c,c])\calH_a,\product_a)$ is a Hilbert space. Since this space has finite codimension in $\calH_a$, we obtain $\kappa_-(\product_a,\calH_a) < \infty$, which is obviously impossible.
\end{proof}

Standard perturbation-theoretical arguments imply that also for non-real $z$ the operators $A(z)$ are closed and densely defined, and for all $\la\in\C$ the operator $A(z) - \la$ is Fredholm with index zero.

The following abbreviations will be used throughout this paper:
$$
\calI := {\bf [-\pi,\pi]},\quad\calH_a := L^2_{|w|}(0,a)\quad\text{and}\quad\wt\calH := L^2(\calI,\calH_a).
$$
The {\it multiplication operator} $\wt T$ with the family of self-adjoint operators $\{T(t) : t\in\calI\}$ is an operator in $\wt\calH$. It has the domain of definition
$$
\dom\wt T := \{F\in\wt\calH : F(t)\in\dom T(t)\text{ for a.e.\ $t\in\calI$ and }T(\cdot)F(\cdot)\in\wt\calH\}
$$
and acts in the following way:
\begin{equation}\label{e:wtT}
\big(\wt T F\big)(t) := T(t)F(t),\quad F\in\dom\wt T,\;t\in\calI.
\end{equation}
The {\it Gelfand transform} $\calG : L^2_{|w|}(\R)\,\to\,\wt\calH$ is defined by
$$
(\calG f)(t) := \liim\limits_{N\to\infty}\;\frac 1 {\sqrt{2\pi}}\sum_{n=-N}^N \!e^{-int}f(\,\cdot+na),\quad t\in\calI,\,f\in L^2_{|w|}(\R).
$$
Here, $\liim$ denotes the limit in $\wt\calH = L^2(\calI,\calH_a)$. It is well-known (see, e.g., \cite[Lemma 16.7 and Satz 16.9]{w2}) that the Gelfand transform is unitary and that
\begin{equation}\label{e:Tmop}
\wt T = \calG T\calG^{-1}.
\end{equation}
In particular, $\wt T$ is a self-adjoint operator in $\wt\calH$. Moreover, the operator $\wt J$ in $\wt\calH$, given by
$$
\big(\wt J F\big)(t) := J_aF(t),\quad F\in\wt\calH,\;t\in\calI,
$$
is a fundamental symmetry in $\wt\calH$. For $f\in L^2_{|w|}(\R)$ and $t\in\calI$ we have
\begin{align*}
\big(\wt J \calG f\big)(t)
&= \liim_{N\to\infty}\,\frac 1 {\sqrt{2\pi}} \sum_{n=-N}^N\,e^{-int}\sgn(w(\cdot))f(\cdot+na)\\
&= \liim_{N\to\infty}\,\frac 1 {\sqrt{2\pi}} \sum_{n=-N}^N\,e^{-int}\sgn(w(\cdot+na))f(\cdot+na)\\
&= \liim_{N\to\infty}\,\frac 1 {\sqrt{2\pi}} \sum_{n=-N}^N\,e^{-int}(Jf)(\cdot+na)\\
&= \left(\calG Jf\right)(t),
\end{align*}
and thus
\begin{equation}\label{e:UJ}
\wt J = \calG J\calG^{-1}.
\end{equation}
Let $\wt A$ be the multiplication operator with the family $\{A(t) : t\in\calI\}$, i.e.
$$
\dom\wt A = \{F\in\wt\calH : F(t)\in\dom A(t)\text{ for a.e.\ $t\in\calI$ and }A(\cdot)F(\cdot)\in\wt\calH\},
$$
and
$$
\big(\wt A F\big)(t) = A(t)F(t),\quad F\in\dom\wt A,\;t\in\calI.
$$
Then $\dom\wt A = \dom\wt T$, and \eqref{e:Tmop}--\eqref{e:UJ} imply
$$
\wt A = \wt J\wt T = \calG JT\calG^{-1} = \calG A\calG^{-1}.
$$
We summarize the above discussion in the following lemma.

\begin{lem}\label{l:uebersicht}
Let $\wt A$, $\wt T$ and $\wt J$ be the multiplication operators in $\wt\calH$ with the families of operators $\{A(t) : t\in\calI\}$, $\{T(t) : t\in\calI\}$ and $\{J_a\}$, respectively. Then with the Gelfand transform $\calG$ the following holds:
$$
\wt A = \calG A\calG^{-1},\quad \wt T = \calG T\calG^{-1}\quad\text{and}\quad\wt J = \calG J\calG^{-1}.
$$
\end{lem}

In the following we introduce the Floquet discriminant and the monodromy matrix of $\fraka$. For $\la\in\C$ denote by $\vphi_\la$ and $\psi_\la$ the solutions of $\fraka(u) = \la u$ which satisfy the initial conditions
\begin{align}
\begin{split}\label{e:phipsi}
\vphi_\la(0) = 1,\quad &(p\vphi_\la')(0) = 0,\\
\psi_\la(0) = 0,\quad &(p\psi_\la')(0) = 1.
\end{split}
\end{align}
We mention that $\vphi_\la(x)$, $\psi_\la(x)$, $p\vphi_\la'(x)$ and $p\psi_\la'(x)$ are entire functions (in $\la\in\C$) for every $x\in\R$ and that $(\la,x)\mapsto (\vphi_\la(x),\psi_\la(x),(p\vphi_\la')(x),(p\psi_\la')(x))$ is continuous on $\C\times\R$. The entire function
$$
D(\la) := \vphi_\la(a) + (p\psi_\la')(a)
$$
is called the {\em Floquet discriminant} of $\fraka$. Here, since
\begin{equation}\label{e:phipsi_quer}
\vphi_{\ol\la} = \ol{\vphi_\la}\quad\text{and}\quad\psi_{\ol\la} = \ol{\psi_\la},
\end{equation}
it has the additional property $D(\ol\la) = \ol{D(\la)}$, $\la\in\C$. In particular, if $\la$ is real then $\vphi_\la$ and $\psi_\la$ are real-valued and $D(\la)$ is real. The Floquet discriminant is the trace of the so-called {\it monodromy matrix}
\begin{equation}\label{e:L}
L(\la) := \begin{pmatrix}\vphi_\la(a) & \psi_\la(a)\\(p\vphi_\la')(a) & (p\psi_\la')(a)\end{pmatrix}.
\end{equation}
As $\det L(\la) = 1$ for all $\la\in\C$ we have
\begin{equation}\label{e:mu_muInv}
\rho\in\sigma(L(\la))\quad\Llra\quad\rho^{-1}\in\sigma(L(\la)).
\end{equation}
The following lemma is well-known, see, e.g., \cite{e}.

\begin{lem}\label{l:DundAt}
Let $\la,z\in\C$. Then,
$$
\la\in\sigma(A(z))\quad\Llra\quad \sigma(L(\la)) = \{e^{iz},e^{-iz}\}\quad\Llra\quad D(\la) = 2\cos(z).
$$
In particular, the spectra of the operators $A(t)$, $t\in [0,\pi]$, are mutually disjoint.
\end{lem}

We point out an important fact which is a consequence of Lemma \ref{l:DundAt} and Proposition \ref{p:cl}.

\begin{cor}\label{c:nc}
The entire function $D$ is not a constant.
\end{cor}

Indeed, if $D$ is a constant, then Lemma \ref{l:DundAt} implies that $\sigma(A(t)) = \emptyset$ for all but at most two $t\in\calI$. If the weight function $w$ is indefinite, this contradicts Proposition \ref{p:cl}. Otherwise, $J_a = \pm I$, and for each $t\in\calI$ the operator $A(t) = \pm T(t)$ is self-adjoint and hence has a non-empty spectrum.

In the next lemma we consider the function $(z,\la)\mapsto R(z,\la)$, where
$$
R(z,\la) := (A(z) - \la)^{-1}.
$$
Note that by Lemma \ref{l:DundAt} the resolvent $R(z,\la)$ exists if and only if $D(\la)\neq 2\cos(z)$.

\begin{lem}\label{l:res}
The following statements hold:
\begin{enumerate}
\item[{\rm (a)}] For each $r > 0$ there exists $\la_0\in\C$ such that $\la_0\in\rho(A(z))$ for all $z\in\R + i[-r,r]$.
\item[{\rm (b)}] For fixed $\la\in\C$ the mapping $z\mapsto R(z,\la)$ is holomorphic on the open set $\{z : 2\cos(z)\neq D(\la)\}$.
\item[{\rm (c)}] If $\calU$ is a domain in $\C$ such that $\calU\subset\rho(A(z))$ for all $z$ in a compact set $K\subset\C$, then $R(z,\la)$ is continuous on $K\times\calU$.
\end{enumerate}
\end{lem}
\begin{proof}
Let $r > 0$ and suppose that $\la_0$ as in (a) does not exist. Then for every $\la\in\C$ there exists $z\in\R + i[-r,r]$ such that $\la\in\sigma(A(z))$, or equivalently, $D(\la) = 2\cos(z)$, cf.\ Lemma \ref{l:DundAt}. Hence, the entire function $D$ is bounded and therefore constant which is impossible due to Corollary \ref{c:nc}.

For the proof of (b) let $\la\in\C$. For $z\in\C$ with $d(z,\la) := 2\cos(z) - D(\la)\neq 0$ the resolvent $R(z,\la)$ of $A(z)$ in $\la$ is given by
\begin{equation}\label{e:res1}
\big(R(z,\la)g\big)(x) = \int_0^a\,G_\la(z,x,y)g(y)w(y)\,dy\,,\quad g\in\calH_a,\;x\in [0,a],
\end{equation}
where
\begin{equation}\label{e:res2}
G_\la(z,x,y) = \Psi_\la(x)^T\left(\frac{L(\la) - e^{-iz}}{d(z,\la)} + {\bf 1}_{\{y\le x\}}(x,y)\right)\frakJ\Psi_\la(y).
\end{equation}
Hereby,
$$
\frakJ := \mat 0 1 {-1} 0\quad\text{and}\quad\Psi_\la(x) := \vek{\vphi_\la(x)}{\psi_\la(x)}.
$$
Let $z,\zeta\in\C$ with $d(z,\la)\neq 0$ and $d(\zeta,\la)\neq 0$. Then
$$
G_\la(z,x,y) - G_\la(\zeta,x,y) = \frac{e^{-iz} - e^{-i\zeta}}{d(z,\la)d(\zeta,\la)}\Psi_\la(x)^TM_\la(z,\zeta)\frakJ\Psi_\la(y),
$$
where $M_\la(z,\zeta) := (e^{i(z+\zeta)} - 1)L(\la) + D(\la) - e^{iz} - e^{i\zeta}$. Hence, for $g\in\calH_a$ we have
$$
R(z,\la)g - R(\zeta,\la)g = \frac{e^{-iz} - e^{-i\zeta}}{d(z,\la)d(\zeta,\la)}\Psi_\la^TM_\la(z,\zeta)\frakJ\vek{[g,\vphi_{\ol\la}]_a}{[g,\psi_{\ol\la}]_a}.
$$
This proves (b). And due to
$$
R(z,\la) - R(z_0,\la_0) = \big(R(z,\la) - R(z_0,\la)\big) + \big(R(z_0,\la) - R(z_0,\la_0)\big)
$$
also (c) is proved.
\end{proof}

\begin{rem}
Lemmas \ref{l:DundAt} and \ref{l:res} imply that the operator function $z\mapsto A(z)$ is holomorphic on $\C$ in the sense of \cite[page 366]{k}, cf.\ \cite[Theorem VII-1.3]{k}.
\end{rem}

In the sequel a continuous mapping $\gamma : \calJ\to\C$, where $\calJ$ is a real (bounded or unbounded) interval, will be called a {\it curve}. As usual, we identify $\gamma$ with its image $\gamma(\calJ)$. We shall call a curve $\gamma$ {\it analytic} if the mapping $\gamma : \calJ\to\C$ is injective and analytic at each $t$ in the real interior of $\calJ$.

\begin{thm}\label{t:A_first}
The operator $A$ has the following spectral properties:
\begin{enumerate}
\item[{\rm (i)}]   $\sigma(A) = \bigcup_{t\in\calI}\,\sigma(A(t)) = \bigcup_{t\in [0,\pi]}\,\sigma(A(t)) = \{\la\in\C : D(\la)\in [-2,2]\}$.
\item[{\rm (ii)}]  $\sigma(A) = \sigma_c(A)$.
\item[{\rm (iii)}] $\sigma(A)$ contains neither interior nor isolated points.
\item[{\rm (iv)}]  $\sigma(A)$ consists of closures of analytic curves.
\item[{\rm (v)}]   $\rho(A)$ does not have bounded connected components.
\item[{\rm (vi)}]  If the weight function $w$ is indefinite, then the real spectrum of $A$ is neither bounded from above nor from below.
\end{enumerate}
\end{thm}
\begin{proof}
By Lemma \ref{l:uebersicht} it suffices to prove the theorem for the multiplication operator $\wt A$ instead of $A$. The second and third equality in (i) follow from Lemma \ref{l:DundAt}. Choose some $\la_0$ as in Lemma \ref{l:res} and note that the multiplication operator with the family $\{(A(t) - \la_0)^{-1} : t\in\calI\}$ coincides with $(\wt A - \la_0)^{-1}$. From \cite{dmt} we conclude that
$$
\sigma\big((\wt A - \la_0)^{-1}\big) = \bigcup_{t\in\calI}\,\sigma\big((A(t) - \la_0)^{-1}\big),
$$
which implies (i). Statement (vi) is a consequence of (i) and Proposition \ref{p:cl}. Also (iv) follows from (i). Interior points of $\sigma(\wt A)$ cannot exist according to Corollary \ref{c:nc}, and the absence of isolated points of $\sigma(\wt A)$ follows from the properties of holomorphic functions. If the operator $\wt A$ has an eigenvalue $\la\in\C$, then there exists $F\in\dom\wt A$, $F\neq 0$, such that $A(t)F(t) = \la F(t)$ for a.e.\ $t\in\calI$. But as $D(\la) = 2\cos(t)$ is only possible for at most two $t\in\calI$, from Lemma \ref{l:DundAt} we obtain the contradiction $F = 0$. Statement (ii) now follows from the implication (cf.\ \eqref{e:ran_ker})
$$
\la\in\sigma_r(\wt A)\quad\Lra\quad\ol\la\in\sigma_p(\wt A).
$$
It remains to prove (v). Assume that there exists a bounded connected component of $\rho(A)$. Then $\Im D = 0$ on its boundary. As $\Im D$ is a harmonic function, it follows from the maximum and minimum principle that $\Im D = 0$ in the whole component which contradicts Corollary \ref{c:nc}.
\end{proof}

The following corollary generalizes the main result of \cite{mz}, where the authors assume that the self-adjoint operator $T = JA$ is uniformly positive (the so-called left-definite case). The complex derivative of a function $f : \calU\to\C$, $\calU\subset\C$, is denoted by $\dot f$.

\begin{cor}\label{c:mz}
The operator $A$ is $J$-nonnegative if and only if all operators $A(t)$, $t\in\calI$, are $J_a$-nonnegative. In this case the spectrum of $A$ is real and consists of compact intervals $[\alpha,\beta]$ {\rm (}which might intersect in their endpoints only{\rm )}, such that $D(\alpha) = \pm 2$, $D(\beta) = \mp 2$ and $\pm\dot D(\la) < 0$ for $\la\in (\alpha,\beta)$.
\end{cor}
\begin{proof}
The operator $A$ is $J$-nonnegative if and only if the self-adjoint operator $T = JA$ in the Hilbert space $L^2_{|w|}(\R)$ is non-negative. By Theorem \ref{t:A_first}(i) (applied to $\frakt$ instead of $\fraka$) this is the case if and only if all self-adjoint operators $T(t) = J_aA(t)$, $t\in\calI$, are non-negative in the Hilbert space $L^2_{|w|}(0,a)$. This proves the first statement. Assume now that $A$ is $J$-nonnegative. Then $\sigma(A)$ is real since $\sigma(A(t))$ is real for any $t\in\calI$, cf.\ Proposition \ref{p:kns}. It remains to show that there are no real points $\la_0$ such that $D(\la_0)\in (-2,2)$ and $\dot D(\la_0) = 0$. Suppose that $\la_0$ is such a point. Then for each $\veps > 0$ sufficiently small both equations $D(\la) = D(\la_0)\pm\veps$ have two solutions close to $\la_0$, respectively, and these must be real. Therefore, $\la_0$ is both a maximum and a minimum of $D|\R$ which contradicts Corollary \ref{c:nc}.
\end{proof}

\section{Non-real spectrum and sign types}\label{s:closer}
In this section it is our aim to prove the following theorem which can be regarded as the main result of this paper.

\begin{thm}\label{t:dt}
The non-real spectrum of the operator $A$ is bounded. Moreover, there exists a finite number of points $\la_1,\ldots,\la_n\in\R$, $\la_{j-1} < \la_j$, $j=2,\ldots,n$, such that the following holds:
\begin{enumerate}
\item[{\rm (i)}]   The interval $(-\infty,\la_1)$ is of negative type with respect to $A$.
\item[{\rm (ii)}]  Each interval $(\la_{j-1},\la_j)$, $j=2,\ldots,n$, is of definite type with respect to $A$.
\item[{\rm (iii)}] The interval $(\la_n,\infty)$ is of positive type with respect to $A$.
\end{enumerate}
\end{thm}

\begin{rem}
We mention that adjacent intervals in Theorem \ref{t:dt} might be of the same sign type with respect to $A$. As the following lemmas will reveal, this happens if $\dot D$ vanishes in the common endpoint $\la_j$ of the intervals and if the function $\la\mapsto\dot D(\la)\psi_\la(a)$ does not change its sign in a neighborhood of $\la_j$.
\end{rem}

The statements of Theorem \ref{t:dt} follow immediately from the next three lemmas.

\begin{lem}\label{l:spp_A(t)}
Let $t\in\calI$ and $\la\in\sigma(A(t))\cap\R$. If $\dot D(\la)\neq 0$, then $\psi_\la(a)\neq 0$ or $(p\vphi_\la')(a)\neq 0$, and the following statements hold:
\begin{enumerate}
\item[{\rm (i)}]  If $\psi_\la(a)\neq 0$ then
$$
\la\in\sigma_\pm(A(t))\quad\Llra\quad \pm\dot D(\la)\psi_\la(a) < 0.
$$
\item[{\rm (ii)}] If $(p\vphi_\la')(a)\neq 0$ then
$$
\la\in\sigma_\pm(A(t))\quad\Llra\quad \pm\dot D(\la)(p\vphi_\la')(a) > 0.
$$
\end{enumerate}
In particular, if only $\dot D(\la)\neq 0$, then $\la$ is a spectral point of definite type of $A(t)$.
\end{lem}

\begin{lem}\label{l:sigma_ex}
There exists $\calR > 0$ such that for all $t\in\calI$ the following holds:
\begin{enumerate}
\item[{\rm (i)}]   The non-real spectrum of $A(t)$ is contained in $B_\calR(0)$.
\item[{\rm (ii)}]  The interval $(-\infty,-\calR)$ is of negative type with respect to $A(t)$.
\item[{\rm (iii)}] The interval $(\calR,\infty)$ is of positive type with respect to $A(t)$.
\end{enumerate}
\end{lem}

\begin{lem}\label{l:spp_uebersetzung}
If $\Delta\subset\R$ is of positive {\rm (}negative{\rm )} type with respect to $A(t)$ for all $t\in\calI$, then $\Delta$ is of positive type {\rm (}negative type, respectively{\rm )} with respect to $A$.
\end{lem}

In the proof of Lemma \ref{l:spp_A(t)} (more precisely, in that of Lemma \ref{l:D_and_f} below) we will make use of the relation
\begin{equation}\label{e:D_Strich}
\dot D(\la) = - \psi_\la(a)\,[\vphi_\la,\ol{\vphi_\la}]_a + \big(\vphi_\la(a) - (p\psi_\la')(a)\big)\,[\vphi_\la,\ol{\psi_\la}]_a + (p\vphi_\la')(a)\,[\psi_\la,\ol{\psi_\la}]_a\,.
\end{equation}
To prove \eqref{e:D_Strich}, note that for $g\in\calH_a = L^2_{|w|}(0,a)$ the solution of the initial value problem
$$
\fraka(u) - \la u = g,\quad u(0) = (pu')(0) = 0,
$$
is given by
$$
u(x) = \vphi_\la(x)\cdot\int_0^x\,\psi_\la g w\,dy \,-\, \psi_\la(x)\cdot\int_0^x\,\vphi_\la g w\,dy.
$$
As for $\la,\mu\in\C$ the function $f := \vphi_\la - \vphi_\mu$ has the properties $f(0) = (pf')(0) = 0$ and $\fraka(f) - \la f = (\la - \mu)\vphi_\mu$, it follows that
$$
\frac{\vphi_\la(x) - \vphi_\mu(x)}{\la - \mu} = \vphi_\la(x)\cdot\int_0^x\,\psi_\la \vphi_\mu w\,dy \,-\, \psi_\la(x)\cdot\int_0^x\,\vphi_\la \vphi_\mu w\,dy.
$$
Analogously, we proceed with the functions $\psi_\la$, $p\vphi_\la'$ and $p\psi_\la'$ and obtain the formulas
\begin{align}
\begin{split}\label{e:Strich_formulas}
\frac d {d\la}\,\vphi_\la(a) &= \vphi_\la(a)\,[\psi_\la,\ol{\vphi_\la}]_a - \psi_\la(a)\,[\vphi_\la,\ol{\vphi_\la}]_a\\
\frac d {d\la}\,(p\vphi_\la')(a) &= (p\vphi_\la')(a)\,[\psi_\la,\ol{\vphi_\la}]_a - (p\psi_\la')(a)\,[\vphi_\la,\ol{\vphi_\la}]_a\\
\frac d {d\la}\,\psi_\la(a) &= \vphi_\la(a)\,[\psi_\la,\ol{\psi_\la}]_a - \psi_\la(a)\,[\vphi_\la,\ol{\psi_\la}]_a\\
\frac d {d\la}\,(p\psi_\la')(a) &= (p\vphi_\la')(a)\,[\psi_\la,\ol{\psi_\la}]_a - (p\psi_\la')(a)\,[\vphi_\la,\ol{\psi_\la}]_a.
\end{split}
\end{align}
These imply \eqref{e:D_Strich}.

\begin{lem}\label{l:D_and_f}
Let $t\in [-\pi,\pi]$ and $\la\in\sigma(A(t))$. Moreover, let $f_\la$ and $f_{\ol\la}$ be eigenfunctions of $A(t)$ corresponding to the eigenvalues $\la$ and $\ol\la$, respectively. Then we have
$$
\psi_\la(a)\,[f_\la,f_{\ol\la}]_a = -f_\la(0)\ol{f_{\ol\la}(0)}\cdot\dot D(\la)
$$
and
$$
(p\vphi_\la')(a)\,[f_\la,f_{\ol\la}]_a = (pf_\la')(0)\ol{(pf_{\ol\la}')(0)}\cdot\dot D(\la).
$$
\end{lem}
\begin{proof}
As $\vphi_{\ol\la} = \ol{\vphi_\la}$ and $\psi_{\ol\la} = \ol{\psi_\la}$, there are $\alpha,\beta,\gamma,\delta\in\C$ such that
$$
f_\la = \alpha\vphi_\la + \beta\psi_\la\quad\text{and}\quad f_{\ol\la} = \gamma\ol{\vphi_\la} + \delta\ol{\psi_\la}.
$$
It is obvious that
\begin{equation}\label{e:werte}
\alpha = f_\la(0), \quad \beta = (pf_\la')(0), \quad \gamma = f_{\ol\la}(0), \quad \delta = (pf_{\ol\la}')(0).
\end{equation}
For simplicity, we set $\vphi := \vphi_\la$ and $\psi := \psi_\la$. From $f_\la,f_{\ol\la}\in\dom A(t)$ we deduce the four equations
\begin{align*}
\left(\vphi(a) - e^{it}\right)\alpha &= -\psi(a)\beta\\
\left((p\psi')(a) - e^{it}\right)\beta &= -(p\vphi')(a)\alpha\\
\left(\vphi(a) - e^{-it}\right)\ol\gamma &= -\psi(a)\ol\delta\\
\left((p\psi')(a) - e^{-it}\right)\ol\delta &= -(p\vphi')(a)\ol\gamma.
\end{align*}
With the help of \eqref{e:D_Strich} we obtain
\begin{align*}
-\psi(a)\,[f_\la,f_{\ol\la}]_a
&= -\psi(a)\,\Big( \alpha\ol\gamma\,[\vphi,\ol\vphi]_a + (\beta\ol\gamma\ + \alpha\ol\delta)\,[\vphi,\ol\psi]_a + \beta\ol\delta\,[\psi,\ol\psi]_a\Big)\\
&= -\alpha\ol\gamma\,\psi(a)\,[\vphi,\ol\vphi]_a + \alpha\ol\gamma\left( 2\vphi(a) - e^{it} - e^{-it} \right)\,[\vphi,\ol\psi]_a\\
&\hspace{5.5cm} + \alpha\left(\vphi(a) - e^{it}\right)\ol\delta\,[\psi,\ol\psi]_a\\
&= -\alpha\ol\gamma\,\psi(a)\,[\vphi,\ol\vphi]_a + \alpha\ol\gamma\big(\vphi(a) - (p\psi')(a)\big)\,[\vphi,\ol\psi]_a\\
&\hspace{4cm} + \alpha\left( D(\la) - (p\psi')(a) - e^{it} \right)\ol\delta\,[\psi,\ol\psi]_a\\
&= \alpha\ol\gamma\,\big(\dot D(\la) - (p\vphi')(a)\,[\psi,\ol\psi]_a\big) + \alpha\hspace*{-0.03cm}\left( e^{-it} - (p\psi')(a) \right)\ol\delta\,[\psi,\ol\psi]_a\\
&= \alpha\ol\gamma\,\dot D(\la)
\end{align*}
as well as
\begin{align*}
(p\vphi')(a)\,[f_\la,f_{\ol\la}]_a
&= (p\vphi')(a)\,\Big( \alpha\ol\gamma\,[\vphi,\ol\vphi]_a + (\beta\ol\gamma\ + \alpha\ol\delta)\,[\vphi,\ol\psi]_a + \beta\ol\delta\,[\psi,\ol\psi]_a\Big)\\
&= \left(e^{it} - (p\psi')(a)\right)\beta\ol\gamma\,[\vphi,\ol\vphi]_a\\
&\hspace{.5cm}+ \Big(\beta\left(e^{-it} - (p\psi')(a)\right)\ol\delta - \left((p\psi')(a) - e^{it}\right)\beta\ol\delta\Big)\,[\vphi,\ol\psi]_a\\
&\hspace{.5cm}+ \beta\ol\delta\,(p\vphi')(a)\,[\psi,\ol\psi]_a\\
&= \left(\vphi(a) - e^{-it}\right)\beta\ol\gamma\,[\vphi,\ol\vphi]_a + \beta\ol\delta\,\big(\vphi(a) - (p\psi')(a)\big)\,[\vphi,\ol\psi]_a\\
&\hspace{4cm}+ \beta\ol\delta\,(p\vphi')(a)\,[\psi,\ol\psi]_a\\
&= \beta\ol\delta\,\dot D(\la).
\end{align*}
The assertion now follows from \eqref{e:werte}.
\end{proof}

We are now ready to prove Lemmas \ref{l:spp_A(t)}--\ref{l:spp_uebersetzung}.

\begin{proof}[Proof of Lemma {\rm\ref{l:spp_A(t)}}]
Suppose that $\dot D(\la)\neq 0$ but $\psi_\la(a) = (p\vphi_\la')(a) = 0$. Then the monodromy matrix $L(\la)$ in \eqref{e:L} is a diagonal matrix and hence has its eigenvalues $e^{it}$ and $e^{-it}$ on the diagonal. Since the functions $\vphi_\la$ and $\psi_\la$ are real-valued, it follows that $\vphi_\la(a) = (p\psi_\la')(a) = \pm 1$. But in view of \eqref{e:D_Strich} this implies $\dot D(\la) = 0$. A contradiction.

Assume, e.g., $\psi_\la(a)\neq 0$ and let $f_\la$ be any eigenfunction of $A(t)$ corresponding to the eigenvalue $\la$. Then $f_\la(0)\neq 0$ since otherwise $f_\la = \nu\psi_\la$ with some $\nu\in\C\setminus\{0\}$ and thus $\psi_\la(a) = e^{it}\psi_\la(0) = 0$. From Lemma \ref{l:D_and_f} it follows that
$$
[f_\la,f_\la]_a = -|f_\la(0)|^2\frac{\dot D(\la)}{\psi_\la(a)}
$$
which proves the equivalence in (i). The statement (ii) is proved similarly.
\end{proof}

For $r > 0$ and $\la\in\C$ by $B_r(\la)$ we denote the open disc in the complex plane with center $\la$ and radius $r$.

\begin{proof}[Proof of Lemma {\rm\ref{l:sigma_ex}}]
First of all note that $\calR > 0$ as in Lemma \ref{l:sigma_ex} exists if and only if the set
$$
\bigcup_{t\in\calI}\,\sigma_{\rm ex}(A(t))
$$
is bounded. Also note that the number of points in each $\sigma_{\rm ex}(A(t))$, $t\in\calI$, cannot exceed $\kappa^*$, cf.\ Proposition \ref{p:cl}. The proof is divided into two steps.

{\bf 1.} In this first step we prove the assertion under the assumption that one of the two following cases holds true:
\begin{enumerate}
\item[(I)]  $D(0)\notin [-2,2]$.
\item[(II)] $D(0)\in (-2,2)$ and $\dot D(0)\neq 0$.
\end{enumerate}
Let (I) or (II) be satisfied. In what follows we show the following claim:
\begin{enumerate}
\item[(C)] For each $t_0\in\calI$ there exists $\delta > 0$ such that for all $t\in (t_0-\delta,t_0+\delta)$ and all $\la\in\sigma_{\rm ex}(A(t))$ we have
\begin{equation}\label{e:dist}
\dist\big(\la,\sigma_{\rm ex}(A(t_0))\big)\le 1.
\end{equation}
\end{enumerate}
Then the assertion follows since $\calI$ is compact.

Let $t_0\in\calI$ be arbitrary and let $\la_1,\ldots,\la_n$ be the exceptional eigenvalues of $A(t_0)$ in $\C^+$ and $\la_{n+1},\ldots,\la_{n+k}$ the {\it non-zero} real exceptional eigenvalues of $A(t_0)$. Choose $\veps\in (0,1)$ such that with $B_j := B_\veps(\la_j)$ the following holds:
\begin{enumerate}
\item[(a)] $\sigma(A(t_0))\cap\ol{B_j} = \{\la_j\}$ for $j=1,\ldots,n+k$,
\item[(b)] $\ol{B_j}\subset\C^+$ for $j=1,\ldots,n$,
\item[(c)] $0\notin\ol{B_j}$ for $j=n+1,\ldots,n+k$.
\item[(d)] $\ol{B_i}\cap\ol{B_j} = \emptyset$ for all $i,j = 1,\ldots,n+k$, $i\neq j$.
\end{enumerate}
Denote by $\Gamma_j$, $j=1,\ldots,n+k$, the boundary of $B_j$. As a consequence of Lemma \ref{l:DundAt} there exists $\delta_1 > 0$ such that $\Gamma_j\subset\rho(A(t))$ for all $t\in\calI_{\delta_1}(t_0) := [t_0 - \delta_1,t_0 + \delta_1]$ and all $j\in\{1,\ldots,n+k\}$. Hence, the Riesz-Dunford projection
$$
P_j(t) := -\frac 1 {2\pi i}\int_{\Gamma_j}(A(t) - \la)^{-1}\,d\la
$$
is well-defined for $j=1,\ldots,n+k$ and $t\in\calI_{\delta_1}(t_0)$. Moreover, according to Lemma \ref{l:res}, each function $P_j$ is continuous on $\calI_{\delta_1}(t_0)$ in the uniform operator topology. Therefore, there exists $\delta\in (0,\delta_1)$ such that
$$
\|P_j(t) - P_j(t_0)\| < 1\quad\text{and}\quad\|P_j(t_0)(P_j(t) - P_j(t_0))P_j(t_0)\| < 1
$$
holds for all $j=1,\ldots,n+k$ and all $t\in\calI_\delta(t_0)$. Hence, \cite[Lemma 2.1]{lt} implies
\begin{equation}\label{e:lt}
\kappa_s(\product_a,P_j(t)\calH_a) = \kappa_s(\product_a,\calL_{\la_j}(A(t_0))),\quad s\in\{+,-,0\},
\end{equation}
for each $j\in\{1,\ldots,n+k\}$ and all $t\in\calI_\delta(t_0)$. Let us now see that for each $j\in\{1,\ldots,n+k\}$ and all $t\in\calI_\delta(t_0)$ this implies
\begin{equation}\label{e:nice}
\sum_{\la\in B_j\setminus\C^-}\!\!\!\!\!
\kappa_-\big([A(t)\cdot,\cdot]_a,\calM_\la(A(t))\big) = \kappa_-\big([A(t_0)\cdot,\cdot]_a,M_{\la_j}(A(t_0))\big).
\end{equation}
If $\la_j\notin\R$, i.e.\ $j\in\{1,\ldots,n\}$, then \eqref{e:nice} follows directly from Lemma \ref{l:signatures}(iv). Let $j > n$ such that $\la_j\in\R^+$. Then Lemma \ref{l:signatures} and \eqref{e:lt} imply
$$
\kappa_-\big([A(t_0)\cdot,\cdot]_a,M_{\la_j}(A(t_0))\big) = \kappa_-(\product_a,\calL_{\la_j}(A(t_0))) = \kappa_-(\product_a,P_j(t)\calH_a)
$$
for all $t\in\calI_\delta(t_0)$. From Lemma \ref{l:signatures}(ii) we obtain
\begin{align*}
\kappa_-(\product_a,P_j(t)\calH_a)
&= \sum_{\la\in B_j\cap\C^+}\!\!\!\!\dim\calL_\la(A(t)) + \sum_{\la\in B_j\cap\R}\!\!\kappa_-(\product_a,\calL_\la(A(t)))\\
&= \sum_{\la\in B_j\setminus\C^-}\!\!\!\!\!\kappa_-\big([A(t)\cdot,\cdot]_a,\calM_\la(A(t))\big).
\end{align*}
Similarly, one proves that \eqref{e:nice} holds for $\la_j\in\R^-$.

If $0\notin\sigma(A(t_0))$, we can choose $\delta > 0$ so small that $0\notin\sigma(A(t))$ for all $t\in\calI_\delta(t_0)$. Then $\kappa(t)$ is constant on $\calI_\delta(t_0)$, and from \eqref{e:nice} it follows that
\begin{align*}
\sum_{j=1}^{n+k}\sum_{\la\in B_j\setminus\C^-}\!\!\!\!\!\kappa_-\big([A(t)\cdot,\cdot]_a,\calM_\la(A(t))\big)
&= \sum_{j=1}^{n+k}\,\kappa_-\big([A(t_0)\cdot,\cdot]_a,M_{\la_j}(A(t_0))\big)\\
&= \kappa(t_0) = \kappa(t),
\end{align*}
which shows that the exceptional eigenvalues of each $A(t)$, $t\in\calI_\delta(t_0)$, are contained in the union of all $B_j$ and $B_j^* := \{\ol\la : \la\in B_j\}$, $j=1,\ldots,n+k$. Therefore, \eqref{e:dist} holds if $0\notin\sigma(A(t_0))$. In particular, the lemma is proved in case (I).

It remains to prove the claim (C) in the case (II) for $t_0\in (0,\pi)$ with $D(0) = 2\cos(t_0)$. The value $-t_0$ needs not to be considered since $f\in\calL_\la(A(t))\,\Llra\,\ol f\in\calL_{\ol\la}(A(-t))$ implies $\sigma_{\rm ex}(A(-t)) = \sigma_{\rm ex}(A(t))$. By Lemma \ref{l:spp_A(t)} either $\psi_0(a)\neq 0$ or $(p\vphi_0')(a)\neq 0$. Without loss of generality we assume $\psi_0(a)\neq 0$. Moreover, Lemma \ref{l:spp_A(t)} implies $0\in\sp(A(t_0))\cup\sm(A(t_0))$. In particular, $0\notin\sigma_{\rm ex}(A(t_0))$, cf.\ Lemma \ref{l:sigma_ex_basic}, and zero is a simple (isolated) eigenvalue of $A(t_0)$, cf.\ \eqref{e:sp_ker}. 

Choose $\veps$ and $\delta$ from above so small that $\calI_\delta(t_0)\subset (0,\pi)$ and such that for $B_0 := B_{\veps}(0)$ the following holds:
\begin{enumerate}
\item[(a')] $\sigma(A(t_0))\cap\ol{B_0} = \{0\}$,
\item[(b')] $\partial B_0\subset\rho(A(t))$ for all $t\in\calI_\delta(t_0)$,
\item[(c')] $\dot D(\la)\psi_\la(a)\neq 0$ for all $\la\in B_0$.
\item[(d')] $\ol{B_0}\cap\ol{B_j} = \emptyset$ for $j=1,\ldots,n+k$,
\end{enumerate}
For $t\in\calI_\delta(t_0)$ let $\la(t)$ be the simple (isolated) eigenvalue of $A(t)$ in $B_0$. Then $\la(t)$ is real since otherwise $\ol{\la(t)}$ is another eigenvalue of $A(t)$ in $B_0$. Moreover, as $D(\la)\psi_\la(a)$ does not change sign on $B_0\cap\R$, we have $\la(t)\in\sigma_\pm(A(t))$ if $0\in\sigma_\pm(A(t_0))$, cf.\ Lemma \ref{l:spp_A(t)}. In addition, from $\dot D(\la(t))\dot\la(t) = -2\sin(t)$ and $\calI_\delta(t_0)\subset (0,\pi)$ we see that $\dot\la(t)\neq 0$ for all $t\in\calI_\delta(t_0)$.

Let $t\in\calI_\delta(t_0)$, $t > t_0$. Then, due to Lemma \ref{l:sigma_ex_basic}, $\la(t)\in\sigma_{\rm ex}(A(t))$ if and only if $\pm\la(t) > 0$ and $\la(t)\in\sigma_\mp(A(t))$. This holds if and only if $\pm\dot\la(t_0) > 0$ and $0\in\sigma_\mp(A(t_0))$. By Lemma \ref{l:spp_A(t)} this is equivalent to $\dot\la(t_0)\dot D(0)\psi_0(a) > 0$. Hence, the relation $\dot D(\la(s))\dot\la(s) = -2\sin(s)$ for $s\in\calI_\delta(t_0)$ yields that $\la(t)\in\sigma_{\rm ex}(A(t))$ for $t > t_0$ if and only if $\psi_0(a) < 0$. By $D_+$ denote the Floquet discriminant corresponding to the differential expression $\frakt$ in \eqref{e:frakt}. Then Lemma \ref{l:spp_A(t)} implies that $\dot D_+(0)\psi_0(a) < 0$. Therefore, $\la(t)\in\sigma_{\rm ex}(A(t))$ for $t > t_0$ if and only if $\dot D_+(0) > 0$. Similarly, one proves that $\la(t)\in\sigma_{\rm ex}(A(t))$ for $t < t_0$ if and only if $\dot D_+(0) < 0$.

Assume that $\dot D_+(0) > 0$. Then for $t\in\calI_\delta(t_0)$ we have
$$
\kappa(t) =
\begin{cases}
\kappa^*     &\text{for }t > t_0\\
\kappa^* - 1 &\text{for }t\le t_0.
\end{cases}
$$
Since $\la(t)\notin\sigma_{\rm ex}(A(t))$ for $t\le t_0$ and $\la(t)\in\sigma_{\rm ex}(A(t))$ for $t > t_0$,
$$
\kappa_-\big([A(t)\cdot,\cdot]_a,\calL_{\la(t)}(A(t))\big) =
\begin{cases}
1 &\text{for }t > t_0\\
0 &\text{for }t\le t_0,
\end{cases}\quad t\in\calI_\delta(t_0).
$$
Hence, \eqref{e:nice} and Proposition \ref{p:kns} imply that for each $t\in\calI_\delta(t_0)$ we have
\begin{align*}
\sum_{j=0}^{n+k}\sum_{\la\in B_j\setminus\C^-}\!\!\!\!\!\kappa_-\big([A(t)\cdot,\cdot]_a,\calM_\la(A(t))\big)
&= \kappa_-\big([A(t)\cdot,\cdot]_a,\calL_{\la(t)}(A(t))\big) + \kappa(t_0)\\
&= \kappa_-\big([A(t)\cdot,\cdot]_a,\calL_{\la(t)}(A(t))\big) + \kappa^* - 1 = \kappa(t).
\end{align*}
Therefore, $\sigma_{\rm ex}(A(t))\subset\bigcup_{j=0}^{n+k}(B_j\cup B_j^*)$. A similar reasoning applies if $\dot D_+(0) < 0$. Hence, the lemma is proved for the cases (I) and (II).

{\bf 2.} Assume now that $0\in\sigma(A)$ and that (II) is not satisfied. Then $\dot D(\veps)\neq 0$ and $D(\veps)\notin\{-2,2\}$ for $\veps > 0$ sufficiently small. By $D_\veps$ denote the Floquet discriminant associated with the (periodic) differential expression
$$
\fraka_\veps(f) := \fraka(f) - \veps f = \frac 1 w\Big((pf')' + (q-\veps w)f\Big).
$$
Then $D_\veps(\la) = D(\la + \veps)$ and thus $\dot D_\veps(0)\neq 0$ and $D_\veps(0)\notin\{-2,2\}$. By the first step of this proof there exists $\calR > 0$ such that $\sigma_{\rm ex}(A(t)-\veps)\subset B_\calR(0)$ for all $t\in\calI$. Hence, for all $t\in\calI$ the non-real spectrum of $A(t) - \veps$ is contained in $B_\calR(0)$, $(\calR,\infty)$ is of positive type with respect to $A(t) - \veps$ and $(-\infty,-\calR)$ is of negative type with respect to $A(t) - \veps$. Consequently, for all $t\in\calI$ the non-real spectrum of $A(t)$ is contained in $B_\calR(\veps)$, $(\calR + \veps,\infty)$ is of positive type with respect to $A(t)$ and $(-\infty,-\calR+\veps)$ is of negative type with respect to $A(t)$. But this means that $\sigma_{\rm ex}(A(t))\subset B_\calR(\veps)$ holds for all $t\in\calI$.
\end{proof}

\begin{proof}[Proof of Lemma {\rm\ref{l:spp_uebersetzung}}]
Let $\la\in\R$ such that $\la\in\sp(A(t))\cup\rho(A(t))$ for all $t\in\calI$. We have to prove that $\la\in\sp(\wt A)\cup\rho(\wt A)$. Then $\la\in\sp(A)\cup\rho(A)$ follows from Lemma \ref{l:uebersicht}. First of all we show that there exists $\veps > 0$ such that for all $t\in\calI$ we have
\begin{equation}\label{e:lmm}
f\in\dom A(t),\quad\|(A(t) - \la)f\|_a\le\veps\|f\|_a\quad\Lra\quad [f,f]_a\ge\veps\|f\|_a^2.
\end{equation}
Suppose that such an $\veps > 0$ does not exist. Then for each $n\in\N$ there exist $t_n\in\calI$ and $f_n\in\dom A(t_n)$ with $\|f_n\|_a = 1$ such that
$$
\|(A(t_n) - \la)f_n\|_a\,\le\,1/n\quad\text{and}\quad [f_n,f_n]_a < 1/n.
$$
It is no restriction to assume that $(t_n)$ converges to some $t\in\calI$. We set
$$
g_n := (A(t) - \la_0)^{-1}(A(t_n) - \la_0)f_n\,\in\,\dom A(t),
$$
where $\la_0\in\rho(A)$ is arbitrary. Due to Lemma \ref{l:res}(b) the expression
$$
g_n - f_n = \big((A(t) - \la_0)^{-1} - (A(t_n) - \la_0)^{-1}\big)(A(t_n) - \la_0)f_n
$$
tends to zero as $n\to\infty$. The same holds for
$$
(A(t) - \la)g_n = (A(t_n) - \la)f_n + (\la_0 - \la)(g_n - f_n).
$$
Therefore, $\la\in\sap(A(t))$ and thus, by assumption, $\la\in\sp(A(t))$, which implies
$$
\liminf_{n\to\infty}\,[f_n,f_n]_a = \liminf_{n\to\infty}\,[g_n,g_n]_a > 0.
$$
But this contradicts $[f_n,f_n]_a < 1/n$. Hence, \eqref{e:lmm} is proved.

Assume that $\la\in\sigma(\wt A)$ (and hence $\la\in\sap(\wt A)$, cf.\ \eqref{e:R_sap}) and let $(F_n)\subset\dom\wt A$ with $\|F_n\|_\sim = 1$ for all $n\in\N$ and $\|(\wt A - \la)F_n\|_\sim\to 0$ as $n\to\infty$, i.e.
$$
a_n := \int_\calI\,\|(A(t) - \la)F_n(t)\|_a^2\,dt\,\to\,0\quad\text{and}\quad\int_\calI\,\|F_n(t)\|_a^2\,dt = 1.
$$
For $n\in\N$ we define the measurable set
$$
M_n := \big\{t\in\calI : \|(A(t) - \la)F_n(t)\|_a \,>\, \veps\|F_n(t)\|_a\big\}.
$$
Then
$$
\int_{M_n}\,\|F_n(t)\|_a^2\,dt\,\le\,\frac 1 {\veps^2}\int_{M_n}\,\|(A(t) - \la)F_n(t)\|_a^2\,dt\,\le\,\frac{a_n}{\veps^2}\,\to\,0
$$
as $n\to\infty$. Moreover, by \eqref{e:lmm},
$$
[F_n,F_n]_\sim = \int_\calI\,[F_n(t),F_n(t)]_a\,dt\,\ge\,\int_{M_n}\,[F_n(t),F_n(t)]_a\,dt + \veps\int_{\calI\setminus M_n}\|F_n(t)\|_a^2\,dt.
$$
And since
$$
\left|\int_{M_n}\,[F_n(t),F_n(t)]_a\,dt\right|\,\le\,\int_{M_n}\,\|F_n(t)\|_a^2\,dt\,\to\,0,
$$
it follows that
$$
\liminf_{n\to\infty}\,[F_n,F_n]_\sim\,\ge\,\veps\,\lim_{n\to\infty}\,\int_{\calI\setminus M_n}\|F_n(t)\|_a^2\,dt = \veps.
$$
This shows $\la\in\sp(\wt A)$.
\end{proof}

The rest of this section is devoted to the study of the Floquet discriminant $D$ on $\R$. Recall that the {\it order} of an entire function $v : \C\to\C$ is defined as the infimum of all $c > 0$ with the property
$$
v(\la) = O\left(e^{|\la|^c}\right)\quad (|\la|\to\infty).
$$
If there exists no such $c > 0$, we say that the function $v$ is of infinite order. A proof of the following lemma can be found in \cite[Section VII.1.1]{ys}.

\begin{lem}\label{l:ordnung}
The order of the entire function $D$ is at most one.
\end{lem}

The next lemma is proved in \cite{c}, see \cite[Lemma XI-3.1]{c}. Note that the additional assumption $f(0) = 1$ in \cite{c} is redundant.

\begin{lem}\label{l:ord_abl}
Let $f : \C\to\C$ be a non-constant entire function whose order is at most one and let the zeros $\la_1,\la_2,\ldots$ of $f$ {\rm (}counting multiplicities{\rm )} be ordered in such a way that $|\la_j|\le |\la_{j+1}|$, $j\in\N$. Then for $\la\notin\{\la_k : k\in\N\}$ we have
$$
\frac{f(\la)\ddot f(\la) - \dot f(\la)^2}{f(\la)^2} = -\sum_{k=1}^\infty\,\frac 1 {(\la - \la_k)^2}.
$$
\end{lem}

In the case of a definite weight function $w$ it is well-known that for real $\la$ with $\dot D(\la) = 0$ we have $|D(\la)|\ge 2$ and $D(\la)\ddot D(\la) < 0$, cf.\ \cite[Theorem 12.7]{w3}. The following proposition shows that in the general case the function $D$ has this behaviour on $\R\setminus (-\calR_0,\calR_0)$, where
$$
\calR_0 := \sqrt 2\,\max\big\{|\la| : \la\in\big(\sigma(A(0))\setminus\R\big)\cup\big(\sigma(A(\pi))\setminus\R\big)\cup\{0\}\big\}.
$$
The constant $\calR_0$ is well-defined due to Proposition \ref{p:cl}.

\begin{prop}\label{p:R}
For each $\la\in\R\setminus (-\calR_0,\calR_0)$ with $\dot D(\la) = 0$ we have
$$
|D(\la)|\ge 2 \quad\text{and}\quad D(\la)\ddot D(\la) < 0.
$$
Consequently, $\la$ is a maximum of $D|\R$ if $D(\la)\ge 2$, and a minimum if $D(\la)\le -2$.
\end{prop}
\begin{proof}
Let $\la\in\R\setminus (-\calR_0,\calR_0)$ such that $\dot D(\la) = 0$. By $\la_1,\ldots,\la_n$ we denote the zeros of the function $D(\cdot) - 2$ (and thus the eigenvalues of $A(0)$, cf.\ Lemma \ref{l:DundAt}) in $\C^+$ and set $\la_{n+j} := \ol{\la_j}$ for $j=1,\ldots,n$. In addition, let $\la_{2n+1},\la_{2n+2},\ldots$ be the (infinitely many) real zeros of $D(\cdot) - 2$ such that $|\la_j|\le |\la_{j+1}|$ for $j\ge 2n+1$. From $|\la|\ge\sqrt{2}\,\max_{j=1,\ldots,n}\,|\la_j|$ it follows that
\begin{align*}
|\la - \Re\,\la_j| &\ge |\la| - |\Re\,\la_j| \ge \sqrt{2((\Re\,\la_j)^2 + (\Im\,\la_j)^2)} - |\Re\,\la_j|\\
                   &\ge \sqrt{(|\Re\,\la_j| + |\Im\,\la_j|)^2} - |\Re\,\la_j| = |\Im\,\la_j|
\end{align*}
for $j=1,\ldots,n$, and with an easy calculation one confirms that this implies
$$
(\la - \la_j)^{-2} + (\la - \ol{\la_j})^{-2}\ge 0
$$
for each $j\in\{1,\ldots,n\}$. We apply Lemma \ref{l:ord_abl} and obtain
$$
D(\la)\neq 2 \quad\Lra\quad (D(\la) - 2)\ddot D(\la) < 0.
$$
An analog treatment of the function $D(\cdot) + 2$ gives
$$
D(\la)\neq -2 \quad\Lra\quad (D(\la) + 2)\ddot D(\la) < 0.
$$
These two implications yield the assertion.
\end{proof}

\begin{cor}\label{c:finite_crossings}
The real accumulation points of the non-real spectrum of $A$ are contained in $(-\calR_0,\calR_0)$ and constitute a finite set.
\end{cor}
\begin{proof}
Let $\la_0\in\R$ be an accumulation point of the non-real spectrum of $A$. As the spectrum of $A$ is symmetric with respect to the real axis, in each neighborhood of $\la_0$ in $\C$ there is a pair $\la,\ol\la\in\C\setminus\R$ such that $D(\la) = D(\ol\la)$ which implies $\dot D(\la_0) = 0$. Hence, in $(-\calR_0,\calR_0)$ there is only a finite number of such accumulation points. Suppose now that $|\la_0|\ge\calR_0$. Then from Proposition \ref{p:R} it follows that $D(\la_0)\in\{-2,2\}$ and $D(\la_0)\ddot D(\la_0) < 0$. If, e.g., $D(\la_0) = 2$, then $D|\R$ has a maximum at $\la_0$ and hence, in each neighborhood of $\la_0$, in addition to $\la$ and $\ol\la$, there is also some $\mu\in\R$ such that $D(\la) = D(\ol\la) = D(\mu)$, which contradicts $\ddot D(\la_0)\neq 0$.
\end{proof}

\begin{cor}\label{c:both_real}
If the spectra of $A(0)$ and $A(\pi)$ are real, then the spectrum of $A$ is real.
\end{cor}
\begin{proof}
If the spectra of $A(0)$ and $A(\pi)$ are real, then $\calR_0 = 0$. Hence, by Corollary \ref{c:finite_crossings} the non-real spectrum of $A$ does not accumulate to any real point. Therefore, since $\sigma(A)\setminus\R$ is bounded by Theorem \ref{t:dt}, the set $K := \sigma(A)\cap\C^+$ is compact. Suppose that $K\neq\emptyset$. Then the function $D$ attains its maximum on $K$ (note that $D$ is real-valued on $K$). Let $\la_0\in K$ such that $D(\la)\le D(\la_0)$ for all $\la\in K$. As $\la_0\notin\R$ we have $D(\la_0)\in (-2,2)$. Let $\calU\subset\C^+$ be an open neighborhood of $\la_0$. Then $D(\calU)$ is an open neighborhood of $D(\la_0)$. Hence, there exists $\la_1\in\calU$ such that $D(\la_0) < D(\la_1) < 2$. Since this also implies $\la_1\in K$, we have obtained a contradiction.
\end{proof}

\section{The spectral function and its singularities}\label{s:sf}
A {\it regular spectral curve} of $A$ is an analytic curve $\gamma : \calJ\to\sigma(A)$ such that the derivative $\dot D$ of $D$ does not vanish on $\gamma$, i.e.\ $\dot D(\gamma(t))\neq 0$ for all $t\in\calJ$. If $\gamma$ is a regular spectral curve of $A$, then for $\Delta\subset\gamma$ we denote by $\calI_\Delta$ the set of all $t\in\calI$ such that $\sigma(A(t))\cap\Delta\neq\emptyset$. Each regular spectral curve of $A$ is bounded due to Proposition \ref{p:cl} and Theorem \ref{t:dt}. We say that a regular spectral curve $\gamma$ of $A$ is {\it maximal}, if for each endpoint $\la$ of $\gamma$ we either have $\dot D(\la) = 0$ (and hence $\la\notin\gamma$) or $\dot D(\la)\neq 0$, $D(\la)\in\{-2,2\}$ and $\la\in\gamma$.

\begin{thm}\label{t:sf}
The operator $A$ has a local spectral function on each of its regular spectral curves.
\end{thm}
\begin{proof}
Let $\gamma$ be a regular spectral curve of $A$. Since each regular spectral curve of $A$ is contained in a maximal one, it is no restriction to assume that $\gamma$ is maximal. Then there exists a closed rectifiable Jordan contour $\Gamma$ such that (see Figure \ref{fig:1})
$$
\sigma(A)\cap\operatorname{int}\Gamma = \gamma
\quad\text{and}\quad
\sigma(A)\cap\Gamma = \ol{\gamma}\setminus\gamma.
$$
\begin{figure*}[ht]
\begin{tikzpicture}[line cap=round,line join=round,x=1.0cm,y=1.0cm]
\clip(-2,1.5) rectangle (6.4,5.5);
\draw [shift={(5.8,-0.6)},line width=1pt]  plot[domain=1.5:2.2,variable=\t]({1*4.59*cos(\t r)+0*4.59*sin(\t r)},{0*4.59*cos(\t r)+1*4.59*sin(\t r)});
\draw [shift={(6.66,4.2)},line width=1pt]  plot[domain=2.55:4.15,variable=\t]({1*1.86*cos(\t r)+0*1.86*sin(\t r)},{0*1.86*cos(\t r)+1*1.86*sin(\t r)});
\draw [rotate around={36.79:(3.35,3.02)},line width=1pt] (3.35,3.02) ellipse (1.73cm and 0.96cm);
\draw (2.4,3.8) node {$\Gamma$};
\draw (4,3.3) node {$\gamma$};
\draw (2.9,2.8) node {$|D| = 2$};
\draw (2.9,2.4) node {$\dot D \neq 0$};
\draw (5.4,4.2) node {$\dot D = 0$};
\pgfcircle[fill]{\pgfxy(3.1,3.1)}{2pt}
\pgfcircle[fill]{\pgfxy(4.83,3.88)}{2pt}
\end{tikzpicture}
\caption{Maximal regular spectral curve $\gamma$ and contour $\Gamma$}\label{fig:1}
\end{figure*}
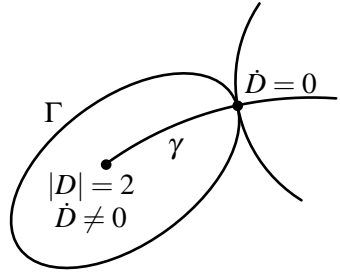

Let $\Delta\in\frakB_0(\gamma)$. Then each of the operators $A(t)$, $t\in\calI_\Delta$, has exactly one eigenvalue in $\operatorname{int}\Gamma$ and there exists a small neighborhood $\calU$ of $\Gamma$ such that $\calU\subset\rho(A(t))$ holds for all $t\in\ol{\calI_\Delta}$. According to Lemma \ref{l:res}(c) the function $R(t,\la)$ is continuous and therefore uniformly bounded on $\ol{\calI_\Delta}\times\Gamma$. Hence, the operator function $t\mapsto E(t;\Delta)$, $t\in\calI$, defined by
\begin{equation}\label{e:e_def}
E(t;\Delta) := 
\begin{cases}
-\frac{1}{2\pi i}\,\int_\Gamma\,R(t,\la)\,d\la\,, &\text{if }t\in\calI_\Delta\\
0, &\text{if }t\in\calI\setminus\calI_\Delta,
\end{cases}
\end{equation}
is measurable and uniformly bounded. Therefore, the multiplication operator $\wt E(\Delta)$ in $\wt\calH$ corresponding to the family $\{E(t;\Delta) : t\in\calI\}$ is an element of $\calB(\wt\calH)$. Since the above construction is independent of the choice of $\Gamma$, the operator $\wt E(\Delta)$ is properly defined. In the following we show that (S1)--(S5) in Definition \ref{d:lsf} hold with $T$, $E$ and $\Xi$ replaced by $\wt A$, $\wt E$ and $\gamma$. Then the theorem is proved according to Lemma \ref{l:uebersicht}.

Let $\Delta\in\frakB_0(\gamma)$. As for each $t\in\calI$ the operator $E(t;\Delta)$ is a projection commuting with the resolvent of $A(t)$, also $\wt E(\Delta)$ is a projection commuting with the resolvent of $\wt A$ (the fact that $\wt E(\Delta)$ is even in the double-commutant of the resolvent of $\wt A$ will be proved below). Moreover, for every $F\in\wt\calH$ the function $\wt E(\Delta)F$ belongs to $\dom\wt A$ since the function $\la R(t,\la)$ is continuous on $\ol{\calI_\Delta}\times\Gamma$. Property (S3) follows directly from the definition of $\wt E$. Let us prove (S2). To this end let $F\in\wt\calH$ and set $\Delta := \bigcup_{j=1}^\infty\Delta_j$ as well as
$$
G_n := \wt E(\Delta)F - \sum_{j=1}^n\,\wt E(\Delta_j)F.
$$
Note that $\calI_\Delta = \bigcup_{j=1}^\infty\calI_{\Delta_j}$ and that the $\calI_{\Delta_j}$ are mutually disjoint. From the definition of $\wt E$ it follows that
$$
G_n(t) =
\begin{cases}
0               & \text{if }t\in\calI\setminus\calI_\Delta,\,n\in\N,\\
E(t;\Delta)F(t) & \text{if }t\in\calI_{\Delta_k},\,n < k,\\
0               & \text{if }t\in\calI_{\Delta_k},\,n\ge k\,.
\end{cases}
$$
Hence, for each $t\in\calI$ we have $G_n(t)\to 0$ as $n\to\infty$. Moreover, as $E(t;\Delta)$ is uniformly bounded in $t\in\calI_\Delta$, there exists $C > 0$ such that $\|G_n(t)\|_a\le C\|F(t)\|_a$ holds for all $t\in\calI$ and all $n\in\N$. Therefore, by Lebesgue's theorem, $\|G_n\|_\sim\to 0$ as $n\to\infty$ as desired.

For the proof of (S4) let $\la_0\in\C\setminus\ol\Delta$ and let $\Gamma_1$ be a rectifiable Jordan contour such that
$$
\ol\Delta\,\subset\,\operatorname{int}\Gamma_1,\quad\la_0\notin\ol{\operatorname{int}\Gamma_1},\quad\sigma(A)\cap\Gamma_1\subset\gamma.
$$
Then $E(t;\Delta)$ is given by \eqref{e:e_def} with $\Gamma$ replaced by $\Gamma_1$. For $G\in\wt E(\Delta)\wt\calH$ set $F(t) = 0$ for $t\in\calI\setminus\calI_\Delta$ and
\begin{equation}\label{e:both}
F(t) := -\frac 1 {2\pi i}\,\int_{\Gamma_1}\,\frac{R(t,\la)G(t)}{\la - \la_0}\,d\la
\end{equation}
for $t\in\calI_\Delta$. Then for a.e.\ $t\in\calI$ we have $F(t)\in E(t;\Delta)\calH_a$ and $\|F(t)\|_a\le C\|G(t)\|_a$ with a constant $C > 0$ which does not depend on $t$. Hence, $F\in\wt E(\Delta)\wt\calH$. Moreover, for a.e.\ $t\in\calI$ we have $(A(t) - \la_0)F(t) = G(t)$ and therefore $F\in\dom\wt A$, $(\wt A - \la_0)F = G$ which shows $\la_0\in\rho(\wt A|\wt E(\Delta)\wt\calH)$.

In order to prove (S5) let $\la_0\notin\ol{\sigma(\wt A)\setminus\Delta}$. Then either $\la_0\in\rho(\wt A)$ or $\la_0$ is in the $\sigma(\wt A)$-interior of $\Delta$. In the first case it is clear that $\la_0\in\rho(\wt A|(I - \wt E(\Delta))\calH)$. In the second case we have $\la_0\in\rho(A(t))$ for all $t\in\ol{\calI\setminus\calI_\Delta}$. Hence, if $G\in (I - \wt E(\Delta))\wt\calH$, then $F(t) := R(t,\la_0)G(t)$ is a proper definition for $t\in\calI\setminus\calI_\Delta$. For $t\in\calI_\Delta$ we define $F(t)$ as in \eqref{e:both} with $\Gamma_1$ replaced by $\Gamma$ (note that $\la_0\in\operatorname{int}\Gamma$). Due to Lemma \ref{l:res}(c) there exists $C > 0$ such that $\|F(t)\|_a\le C\|G(t)\|_a$ for a.e.\ $t\in\calI$. Moreover, $F(t)\in (I - E(t;\Delta))\calH_a$ for a.e.\ $t\in\calI$ implies $F\in (I - \wt E(\Delta))\wt\calH$. In addition, $F(t)\in\dom A(t)$ and $(A(t) - \la_0)F(t) = G(t)$ holds for a.e.\ $t\in\calI$. Consequently, $F\in\dom\wt A$ and $(\wt A - \la_0)F = G$. Therefore, $\la_0\in\rho(\wt A|(I - \wt E(\Delta))\wt\calH)$.

It remains to prove that $\wt E(\Delta)$ is in the double-commutant of the resolvent of $\wt A$ for $\Delta\in\frakB_0(\gamma)$. For this it suffices to consider only closed $\Delta\in\frakB_0(\gamma)$. For $k$ sufficiently large, say $k\ge K$, the set $\Delta_k := \{\la\in\gamma : \dist(\la,\Delta)\le 1/k\}$ is an element of $\frakB_0(\gamma)$. If the spectral curve $\gamma$ is a spectral set of $\wt A$ and $\Delta = \gamma$, then $\wt E(\Delta)$ coincides with the Riesz-Dunford spectral projection of $\wt A$ corresponding to $\Delta$, and hence (S1) holds true. Otherwise, $\Delta$ is a proper subset of $\Delta_k$ for all $k\ge K$. Let $\wt B\in\calB(\wt\calH)$ be an operator which commutes with $\wt S := (\wt A - \la_0)^{-1}$ for some $\la_0\in\rho(\wt A)$. Set $\wt B_0 := \wt B|\wt E(\Delta)\wt\calH\in\calB(\wt E(\Delta)\wt\calH,\wt\calH)$. Then, since $\wt S$ commutes with all $\wt E(\Delta_k)$, we have for $k\ge K$:
$$
\left[\wt S|(I - \wt E(\Delta_k))\wt\calH\right]\,\left[(I - \wt E(\Delta_k))\wt B_0\right] = \left[(I - \wt E(\Delta_k))\wt B_0\right]\,\left[\wt S|\wt E(\Delta)\wt\calH\right].
$$
Owing to (S4), (S5) and Rosenblum's corollary (see, e.g., \cite{rr}) it follows that $(I - \wt E(\Delta_k))\wt B_0 = 0$ for all $k\ge K$ or, equivalently,
$$
\wt B\wt E(\Delta)\wt\calH\subset\wt E(\Delta_k)\wt\calH.
$$
Similarly one proves that for all $k\ge K$
$$
\wt B(I - \wt E(\Delta_k))\wt\calH\subset(I - \wt E(\Delta))\wt\calH.
$$
We will now prove that the following two inclusions hold:
\begin{equation}\label{e:last}
\bigcap_{k\ge K}\wt E(\Delta_k)\wt\calH\,\subset\,\wt E(\Delta)\wt\calH,
\quad
(I - \wt E(\Delta))\wt\calH\,\subset\,\operatorname{c.l.s.}\left\{(I - \wt E(\Delta_k))\wt\calH : k\ge K\right\}.
\end{equation}
Then the proof of (S1) is complete. For simplicity we assume $K = 1$ and set
$$
\delta_0 := \Delta,\;\;\delta_k := \Delta_k\setminus\Delta_{k+1}\;\text{ for }k\ge 1.
$$
Then the $\delta_k$, $k\ge 0$, are mutually disjoint and their union coincides with $\Delta_1$. By (S2) we obtain for every $F\in\wt\calH$:
$$
\wt E(\Delta_1)F = \wt E(\Delta)F + \sum_{k=1}^\infty\,\wt E(\delta_k)F.
$$
As for $k\ge 1$ we have $\wt E(\delta_k) = \wt E(\Delta_k) - \wt E(\Delta_{k+1})$, this implies $\|\wt E(\Delta_k)F - \wt E(\Delta)F\|_\sim\to 0$ as $k\to\infty$ and thus \eqref{e:last}. The stated uniqueness of $\wt E$ is a consequence of (S1), (S2), (S4) and (S5), see, e.g.\ \cite[Lemma 3.14]{j03}.
\end{proof}

Due to Theorem \ref{t:A_first} the spectrum of $A$ consists of the union of countably many regular spectral curves of $A$ and the (discrete) set $c(A)$ of points $\la\in\sigma(A)$ for which $\dot D(\la) = 0$. The points in $c(A)$ will be called the {\em critical points} of $A$.

\begin{rem}\label{r:gsf}
We mention that in a similar manner to the proof of Theorem \ref{t:sf} a spectral projection $E(\Delta)$ can be defined for small connected $\sigma(A)$-neighborhoods $\Delta$ of the critical points of $A$. Hence, if $\frakB_{c(A)}(\sigma(A))$ denotes the collection of all bounded Borel sets in $\sigma(A)$ whose $\sigma(A)$-boundary does not contain any critical point of $A$, then there exists an operator-valued mapping $E$ on $\frakB_{c(A)}(\sigma(A))$ with the properties (S1)--(S5) in Definition \ref{d:lsf} (with $T$ and $\frakB_0(\gamma)$ replaced by $A$ and $\frakB_{c(A)}(\sigma(A))$, respectively).
\end{rem}

\begin{defn}
The mapping $E$ on $\frakB_{c(A)}(\sigma(A))$ from Remark \ref{r:gsf} is called the {\em spectral function} of $A$ {\rm (}with the set of critical points $c(A)${\rm )}. A critical point $\la_0$ of $A$ is called {\em regular} if for some domain $\calU$ in $\C$ with $\la_0\in\calU$ and $\ol\calU\cap c(A) = \{\la_0\}$ we have
$$
\sup_\gamma\sup_{\Delta\subset\calU,\,\Delta\in\frakB_0(\gamma)}\|E(\Delta)\| < \infty,
$$
where the first supremum runs over all regular spectral curves $\gamma$ of $A$ with $\la_0$ in the $\sigma(A)$-boundary of $\gamma$. If the critical point $\la_0$ is not regular, it is called {\em singular}. A singular critical point of $A$ is also called a {\it singularity} of the spectral function $E$ or a {\it spectral singularity} of $A$.
\end{defn}

If $\gamma$ is a regular spectral curve of $A$, then by $\la_\gamma : \calI_\gamma\to\gamma$ we denote the mapping with $\la_\gamma(t)\in\sigma(A(t))$, $t\in\calI_\gamma$. This mapping is unique and real-analytic: if $\calU$ is a domain in $\C$ with $\calU\cap\sigma(A) = \gamma$ and on which $\dot D$ does not vanish, then for $t\in\calI_\gamma$ we have $\la_\gamma(t) = (D|\calU)^{-1}(2\cos(t))$. By $\wt E$ we denote the spectral function of $\wt A$, i.e.\ $\wt E(\Delta) := \calG E(\Delta)\calG^{-1}$, $\Delta\in\frakB_{c(A)}(\sigma(A))$.

\begin{lem}\label{l:sf_rep}
Let $\gamma$ be a regular spectral curve of $A$ and let $\Delta\in\frakB_0(\gamma)$. Then for all $G,H\in\wt\calH$ we have
$$
\big[\wt E(\Delta)G,H\big]_\sim = \int_{\calI_\Delta}\,
\vek{[\vphi_{\la(t)},H(t)]_a}{[\psi_{\la(t)},H(t)]_a}^T
\frac{L(\la(t)) - e^{-it}}{\dot D(\la(t))}
\vek{[G(t),\ol{\psi_{\la(t)}}]_a}{-[G(t),\ol{\vphi_{\la(t)}}]_a}\,dt,
$$
where $\la(\cdot) := \la_\gamma(\cdot)$ and $L(\cdot)$ is the monodromy matrix from {\rm\eqref{e:L}}.
\end{lem}
\begin{proof}
It is no restriction to assume that $\gamma$ is maximal. By the definition of $\wt E(\Delta)$ in the proof of Theorem \ref{t:sf} we have
$$
\big[\wt E(\Delta)G,H\big]_\sim = -\frac 1 {2\pi i}\,\int_{\calI_\Delta}\int_\Gamma\left[R(t,\la)G(t),H(t)\right]_a\,d\la\,dt,
$$
where $\Gamma$ is a closed rectifiable Jordan contour such that $\sigma(A)\cap\operatorname{int}\Gamma = \gamma$ and $\sigma(A)\cap\Gamma = \ol{\gamma}\setminus\gamma$, cf.\ Figure \ref{fig:1}. We will now make use of the representation \eqref{e:res1}--\eqref{e:res2} of $R(t,\la)$. For $x,y\in [0,a]$ and $t\in\calI$ set $f(t,x,y) := (G(t))(y)w(y)\ol{(H(t))(x)}w(x)$. Since the function
$$
\la\mapsto\int_0^a\int_0^x\,\Psi_\la(x)^T\frakJ\Psi_\la(y)f(t,x,y)\,dy\,dx
$$
is entire for every $t$, it follows that $\big[\wt E(\Delta)G,H\big]_\sim$ coincides with
\begin{align*}
&-\frac 1 {2\pi i}\int_{\calI_\Delta}\int_\Gamma\int_0^a\int_0^a\,
\Psi_\la(x)^T\frac{L(\la) - e^{-it}}{2\cos(t) - D(\la)}\frakJ\Psi_\la(y)f(t,x,y)\,dy\,dx\,d\la\,dt\\
= &-\frac 1 {2\pi i}\int_{\calI_\Delta}\int_\Gamma\,\vek{[\vphi_\la,H(t)]_a}{[\psi_\la,H(t)]_a}^T\frac{L(\la) - e^{-it}}{2\cos(t) - D(\la)}\vek{[G(t),\ol{\psi_\la}]_a}{-[G(t),\ol{\vphi_\la}]_a}\,d\la\,dt.
\end{align*}
The assertion is now a consequence of
$$
\frac{L(\la) - e^{-it}}{2\cos(t) - D(\la)} = -\frac{\frac{\la(t) - \la}{D(\la(t)) - D(\la)}(L(\la) - e^{-it})}{\la - \la(t)}
$$
and Cauchy's integral formula.
\end{proof}

\begin{lem}\label{l:sing}
Let $\la_0\in c(A)$, let $\gamma$ be a regular spectral curve of $A$ with $\ol\gamma\cap c(A) = \{\la_0\}$ and let $f(t)$ be one of the entries of the matrix function $L(\la_\gamma(t)) - e^{-it}$. If
$$
\frac{f}{\dot D\circ\la_\gamma}\notin L^2(\calI_{\gamma}),
$$
then $\la_0$ is a singular critical point of $A$.
\end{lem}
\begin{proof}
Set $\la(\cdot) := \la_\gamma(\cdot)$. We prove the lemma for the entry $f(t) = \vphi_{\la(t)}(a) - e^{-it}$. The proof for the other ones is similar. For linearly independent functions $g,h\in\calH_a$ we set
$$
\Phi(g,h) := J_a\left(g - \frac{(g,h)_a}{\|h\|_a^2}h\right).
$$
Since $f$ is real-analytic, the zeros of $f$ in $\calI_\gamma$ are at most countable. If $f/(\dot D\circ\la)\in L^1(\calI_\gamma)\setminus L^2(\calI_\gamma)$, for $t\in\calI_\gamma$ define
\begin{align*}
G(t) &:= \left|\frac{f(t)}{\dot D(\la(t))}\right|^{1/2}\Phi(\ol{\psi_{\la(t)}},\ol{\vphi_{\la(t)}}),\\
H(t) &:= \frac{f(t)}{\dot D(\la(t))}\left|\frac{f(t)}{\dot D(\la(t))}\right|^{-1/2}\Phi(\vphi_{\la(t)},\psi_{\la(t)}).
\end{align*}
For $t\in\calI\setminus\calI_\gamma$ we set $G(t) := H(t) := 0$. Then, $G,H\in L^2(\calI,\calH_a)$, and due to Lemma \ref{l:sf_rep} for each connected $\Delta\in\frakB_0(\gamma)$ we have
$$
[\wt E(\Delta)G,H]_\sim = \int_{\calI_\Delta}\,\frac{\big(\|\vphi_{\la(t)}\|_a^2\|\psi_{\la(t)}\|_a^2 - |(\vphi_{\la(t)},\psi_{\la(t)})_a|^2\big)^2}{\|\vphi_{\la(t)}\|_a^2\|\psi_{\la(t)}\|_a^2}\cdot\left|\frac{f(t)}{\dot D(\la(t))}\right|^2\,dt\,.
$$
This shows that $[\wt E(\Delta)G,H]_\sim$ tends to $\infty$ when the $\sigma(A)$-boundary of $\Delta$ tends to $\la_0$. If $f/(\dot D\circ\la)\notin L^1(\calI_\gamma)$, the same holds for $[\wt E(\Delta)G,H]_\sim$, where
$$
G(t) := \Phi(\ol{\psi_{\la(t)}},\ol{\vphi_{\la(t)}}),\quad H(t) := \frac{f(t)}{|f(t)|}\cdot\frac{|\dot D(\la(t))|}{\dot D(\la(t))}\Phi(\vphi_{\la(t)},\psi_{\la(t)}),\quad t\in\calI_\gamma.
$$
Hence, in both cases $\la_0$ is a singular critical point of $A$.
\end{proof}

\begin{thm}\label{t:reg}
Let $\la_0$ be a critical point of $A$ and set $t_0 := \arccos(D(\la_0)/2)$. Then the following assertions are equivalent.
\begin{enumerate}
\item[{\rm (i)}]   $\la_0$ is a regular critical point of $A$.
\item[{\rm (ii)}]  $D(\la_0)\in\{-2,2\}$, $\psi_{\la_0}(a) = (p\vphi_{\la_0}')(a) = 0$ and $\ddot D(\la_0)\neq 0$.
\item[{\rm (iii)}] $\ker\big((A(t_0) - \la_0)^2\big) = \ker(A(t_0) - \la_0)$.
\end{enumerate}
\end{thm}
\begin{proof}
Let $m\ge 2$ be the order of $\la_0$ as a zero of the function $D - D(\la_0)$. Then there exists an entire function $F$ with $F(\la_0)\neq 0$ such that $D(\la) - D(\la_0) = (\la - \la_0)^mF(\la)$ for all $\la\in\C$. Hence, we have $(D(\la) - D(\la_0))^{m-1} = (\la - \la_0)^{m(m-1)}F(\la)^m$ and $\dot D(\la)^m = (\la - \la_0)^{m(m-1)}(mF(\la) + (\la - \la_0)\dot F(\la))^m$ for all $\la\in\C$. Combining these identities we obtain
$$
\dot D(\la)^m = C(\la)^m\big(D(\la) - D(\la_0)\big)^{m-1},\quad\text{where }\,C(\la) := \frac{mF(\la) + (\la - \la_0)\dot F(\la)}{F(\la)}\,.
$$
Note that $\lim_{\la\to\la_0}C(\la) = m$. In what follows let $\gamma$ be a regular spectral curve of $A$ with $\ol\gamma\cap c(A) = \{\la_0\}$ and set $\la(\cdot) = \la_\gamma(\cdot)$.

Assume that (i) holds. Then we have
\begin{equation}\label{e:dotD_int}
\big|\dot D(\la(t))\big| = |C(\la(t))|\cdot\big|2(\cos(t) - \cos(t_0))\big|^{\frac{m-1}{m}},\quad t\in\calI_\gamma\,.
\end{equation}
If $D(\la_0)\notin\{-2,2\}$, then either $\vphi_{\la_0}(a)\neq e^{-it_0}$ or $(p\psi_{\la_0}')(a)\neq e^{-it_0}$. Without loss of generality we assume $\vphi_{\la_0}(a)\neq e^{-it_0}$ and $f(t) := \vphi_{\la(t)}(a) - e^{-it}\neq 0$ for $t\in\calI_\gamma$. Then \eqref{e:dotD_int} implies
$$
\frac{f}{\dot D\circ\la}\notin L^2(\calI_\gamma),
$$
which, due to Lemma \ref{l:sing}, is a contradiction. Hence, $D(\la_0)\in\{-2,2\}$. In the following we only consider the case $D(\la_0) = 2$ (and thus $t_0 = 0$). Similar arguments apply to the case $D(\la_0) = -2$. From \eqref{e:dotD_int} it follows that
\begin{equation}\label{e:DD_2}
|\dot D(\la(t))| = c(t)\cdot|t|^{2-2/m},\quad t\in\calI_\gamma,
\end{equation}
where $c\in C(\ol{\calI_\gamma})$ with $c(0)\neq 0$. Hence, if $\psi_{\la_0}(a)\neq 0$ or $(p\vphi_{\la_0}')(a)\neq 0$, then
$$
\frac{\psi_{\la(\cdot)}(a)}{\dot D\circ\la}\notin L^1(\calI_\gamma)
\quad\text{or}\quad
\frac{(p\vphi_{\la(\cdot)}')(a)}{\dot D\circ\la}\notin L^1(\calI_\gamma),
$$
which again contradicts (i). Assume now that $\psi_{\la_0}(a) = (p\vphi_{\la_0}')(a) = 0$, but $\ddot D(\la_0) = 0$. Then $\vphi_{\la_0}(a) = (p\psi_{\la_0}')(a) = 1$ and $m\ge 3$. From $\dot D(\la) = F(\la)C(\la)(\la - \la_0)^{m-1}$ and \eqref{e:DD_2} we obtain
\begin{equation}\label{e:la_2}
|\la(t) - \la_0| = \hat c(t)|t|^{2/m},\quad t\in\calI_\gamma
\end{equation}
with $\hat c\in C(\ol{\calI_\gamma})$, $\hat c(0)\neq 0$. Let $g(\cdot)$ be an entry of $L(\cdot) - I$. We set $f_g(t) := g(\la(t))$ if $g$ is an off-diagonal entry and $f_g(t) := g(\la(t)) + 1 - e^{-it}$ otherwise, $t\in\calI_\gamma$. Then $f_g(t)$ is an entry of $L(\la(t)) - e^{-it}$. By $\kappa(g)$ denote the order of $\la_0$ as a zero of $g$. If $\kappa(g)\le (m-2)/2$, then it is seen from \eqref{e:DD_2} and \eqref{e:la_2} that $f_g/(\dot D\circ\la)\notin L^1(\calI_\gamma)$. Therefore, we have $\kappa(g)> (m-2)/2$ for all entries $g$ of $L(\cdot) - I$. The relation
\begin{equation}\label{e:gl}
-(\vphi_\la(a) - 1)((p\psi_\la')(a) - 1) = D(\la) - 2 + (p\vphi_\la')(a)\psi_\la(a)
\end{equation}
implies that there exists an entry $g$ of $L(\cdot) - I$ with $\kappa(g)\le m/2$. Let $g$ be such an entry of $L(\cdot) - I$. Let us first assume that $m$ is odd. Then $\kappa(g) = (m-1)/2$ and for all $p\in (m/(m-1),m/(m-2))$, $p\le 2$, we have
$$
\frac{g\circ\la}{\dot D\circ\la}\in L^1(\calI_\gamma)\setminus L^p(\calI_\gamma)
\quad\text{and}\quad
\frac{1 - e^{-it}}{\dot D\circ\la}\in L^p(\calI_\gamma),
$$
which implies $f_g/(\dot D\circ\la)\notin L^p(\calI_\gamma)$ and thus $f_g/(\dot D\circ\la)\notin L^2(\calI_\gamma)$, contradicting (i). Therefore, $m$ must be even. For all entries $g$ we have $\kappa(g) > (m-2)/2 = m/2 - 1$ and thus $\kappa(g)\ge m/2$. For off-diagonal entries $g$ even $\kappa(g)\ge m/2 + 1$ holds since otherwise $f_g/(\dot D\circ\la)\notin L^2(\calI_\gamma)$. From this and \eqref{e:gl} it follows that
$$
\kappa(\vphi_\la(a) - 1) = \kappa((p\psi_\la')(a) - 1) = \frac m 2\,.
$$
Set $g_{11}(\la) := \vphi_\la(a) - 1$ and $g_{22}(\la) := (p\psi_\la')(a) - 1$, $\la\in\C$, and $f_{jj} := f_{g_{jj}}$, $j=1,2$. There exist entire functions $c_{jj}$ such that $g_{jj}(\la) = (\la - \la_0)^{m/2}c_{jj}(\la)$ for $\la\in\C$ and $c_{jj}(\la_0)\neq 0$. Note that $c_{11}(\la_0) + c_{22}(\la_0) = D^{(m/2)}(\la_0)/(m/2)! = 0$. Therefore, there exists $j\in\{1,2\}$ such that for $t$ sufficiently close to zero
$$
\left|\frac{(\la(t) - \la_0)^{m/2}c_{jj}(\la(t))- i\sin(t)}{t}\right|\,\ge\,\delta
$$
with some $\delta > 0$. For this $j$ we have
$$
\frac{f_{jj}(t)}{\dot D(\la(t))} = \frac{(\la(t) - \la_0)^{m/2}c_{jj}(\la(t))- i\sin(t)}{t}\cdot\frac{t}{\dot D(\la(t))} + \frac{1 - \cos(t)}{\dot D(\la(t))}.
$$
As
$$
\frac{1 - \cos(t)}{\dot D(\la(t))}\in L^2(\calI_\gamma)
\quad\text{and}\quad
\frac{t}{\dot D(\la(t))}\notin L^2(\calI_\gamma)
$$
it follows that $f_{jj}/(\dot D\circ\la)\notin L^2(\calI_\gamma)$. This finally shows that (i) implies (ii).

Assume that (ii) holds. As above, we only consider the case $D(\la_0) = 2$. From $m = 2$, \eqref{e:DD_2} and \eqref{e:la_2} it follows that each entry of
$$
\frac{L(\la(t)) - e^{-it}}{\dot D(\la(t))}
$$
is bounded as a function of $t\in\calI_\gamma$. Hence, the uniform boundedness of $E(\Delta)$ for $\Delta\in\frakB_0(\gamma)$ is a consequence of Lemma \ref{l:sf_rep}, and (i) follows. Moreover, it is seen from the representation \eqref{e:res1}--\eqref{e:res2} of $R(0,\la) = (A(0) - \la)^{-1}$ that $\la_0$ is a pole of order one of $R(0,\la)$. This yields (iii). Conversely, assume that (iii) is satisfied. Then the spectral subspace of $A(t_0)$ corresponding to the isolated eigenvalue $\la_0$ coincides with $ker(A(t_0) - \la_0)$ and has at most dimension $2$. As due to $\dot D(\la_0) = 0$ the function $D$ attains its values at least twice in a neighborhood $\calU$ of $\la_0$, it follows from \cite[Theorem VII-1.7]{k} that for $t$ close to $t_0$ the operator $A(t)$ has exactly two simple eigenvalues in $\calU$ and hence also $\dim\ker(A(t_0) - \la_0) = 2$. This implies (ii).
\end{proof}

\begin{rem}
We remark that the techniques and results above also apply to differential expressions $\fraka$ with complex-valued coefficients $w$, $p$, $q$ and associated non-constant Floquet-discriminant\footnote{If $w$ is not real-valued, then $\product$, $\product_a$ and $\product_\sim$ are only bounded sesquilinear forms.}. In fact, even for higher order differential expressions with complex-valued coefficients (but without weight) the characterization (iii) in Theorem \ref{t:reg} of (finite) regular critical points was proved by Veliev in \cite{v}. In addition, we mention that in the paper \cite{gt} by F.\ Gesztesy and V.\ Tkachenko also necessary and sufficient conditions have been proved for the point $\infty$ not to be a spectral singularity of a Hill operator (i.e.\ $w = p = 1$) with complex-valued potential $q$. In the proof the authors make use of the asymptotic behaviour of the eigenvalues of the operators $A(t)$ and the functions $\vphi_\la$ and $\psi_\la$. To the best of our knowledge such asymptotics do not exist yet in the case of an indefinite weight function. Therefore, we restrict ourselves to proving only a sufficient condition in the next section.
\end{rem}

\begin{cor}\label{c:sing_finite}
The set of singular critical points of $A$ is finite.
\end{cor}
\begin{proof}
Let $\calR > 0$ be as in Lemma \ref{l:sigma_ex} and let $\la$ be a critical point of $A$ in $\C\setminus B_\calR(0)$. Then $\la\in\sp(A(t))\cup\sm(A(t))$, where $t := \arccos(D(\la)/2)$. Hence, $\ker((A(t) - \la)^2) = \ker(A(t) - \la)$ holds by \eqref{e:sp_ker} which implies that $\la$ is a regular critical point of $A$. Therefore, the singular critical points of $A$ are contained in $B_\calR(0)$. And as any critical point of $A$ is a zero of the non-constant holomorphic function $\dot D$, the statement is proved.
\end{proof}

The following corollary can be found as Theorem 3.8 in \cite{ko}.

\begin{cor}\label{c:nn_cp}
Assume that $A$ is $J$-nonnegative. Then the zero point is the only possible singular critical point of $A$. If $0$ is a critical point of $A$, then it is singular.
\end{cor}
\begin{proof}
Due to Corollary \ref{c:mz} the spectrum of $A$ is real. Since $\sigma_{\rm ex}(A(t)) = \emptyset$ for each $t\in\calI$, the value $\calR > 0$ in Lemma \ref{l:sigma_ex} can be chosen arbitrarily small. Hence, the same arguments as in the proof of Corollary \ref{c:sing_finite} imply that the critical points of $A$ in $\R\setminus\{0\}$ are regular. And since $\psi_0(a) > 0$ (see the proof of \cite[Theorem 12.7]{w3}), the origin is a spectral singularity of $A$ if $0\in c(A)$.
\end{proof}

\section{Regularity of the point \texorpdfstring{$\infty$}{infinity}}\label{s:infty}
We say that $\infty$ is a spectral singularity of $A$ if
$$
\sup_{C > \calR}\|E([\calR,C])\| = \infty\qquad\text{or}\qquad\sup_{C > \calR}\|E([-C,-\calR])\| = \infty,
$$
where $\calR$ is as in Lemma \ref{l:sigma_ex} and $E$ denotes the spectral function of $A$. In this section we provide a condition which ensures that $\infty$ is not a spectral singularity of $A$. A point $x_0\in\R$ is called {\it turning point} of the function $w$ if $w$ is indefinite on $(x_0 - \delta,x_0 + \delta)$ for each $\delta > 0$. In Theorem \ref{t:infty} below we assume that the function $w$ has only finitely many turning points in $[0,a]$ at which it is $1$-simple in the following sense.

\begin{defn}\label{d:1simple}
The function $w$ is called $1${\em -simple} at a turning point $x_0$, if there exist $\delta > 0$, $\tau_+,\tau_- > -1$ and functions $\rho_+\in C^1([x_0,x_0+\delta])$ and $\rho_-\in C^1([x_0 - \delta,x_0])$ with $\rho_+(x_0)\neq 0$, $\rho_-(x_0)\neq 0$ and $\sgn(\rho_+(x_0+x)) = -\sgn(\rho_-(x_0-x)) = \operatorname{const}$ for $x\in [0,\delta]$, such that
$$
w(x) = \rho_\pm(x)|x - x_0|^{\tau_\pm},\quad\pm(x - x_0)\in (0,\delta).
$$
\end{defn}

\begin{rem}
The term "$n$-simple" originates from the paper {\rm\cite{cl}} where ordinary differential expressions of order $2n$, $n\in\N$, were investigated.
\end{rem}
%

\begin{thm}\label{t:infty}
Assume that the function $w$ has only finitely many turning points in $[0,a]$ and that $w$ is $1$-simple at each of them. If $p$ and $p^{-1}$ are essentially bounded in neighborhoods of these turning points, then $\infty$ is not a spectral singularity of $A$.
\end{thm}
\begin{proof}
The proof is divided into two steps.

{\bf 1.} In this step we assume that the self-adjoint operator $T = JA$ is uniformly positive, i.e., there exists $\delta > 0$ such that $(Tf,f)\ge\delta\|f\|^2$ for all $f\in\dom T$. Then also $(T(t)f,f)_a\ge\delta\|f\|_a^2$ for all $t\in\calI$ and all $f\in\dom T(t)$. By $T_{\min}$ denote the minimal operator associated with $\frakt$ on $[0,a]$; that is, $T_{\min}f := \frakt(f)$, $f\in\dom T_{\min}$, where
\begin{align*}
\dom T_{\min} = \{f\in L^2_{|w|}(0,a) : \,
&f,f'\in\AC([0,a]),\\
&f(0) = f(a) = (pf')(0) = (pf')(a) = 0\}.
\end{align*}
Clearly, the symmetric operator $T_{\min}$ is uniformly positive. Denote the Friedrichs- and the Krein-von Neumann extension of $T_{\min}$ by $T_F$ and $T_N$, respectively, and define the sets
\begin{align*}
\calD_N  &:= \{f\in\AC([0,a]) : |f'|^2p\in L^1(0,a)\},\\
\calD(t) &:= \{f\in\calD_N : f(a) = e^{it}f(0)\},\;t\in\calI,\\
\calD_F  &:= \{f\in\calD_N : f(0) = f(a) = 0\}.
\end{align*}
Obviously, $\calD_F\,\subset\,\calD(t)\,\subset\,\calD_N$ for all $t\in\calI$. It is well-known (cf.\ \cite{cl}) that
$$
\dom T_N^{1/2} = \calD_N,\quad\dom T(t)^{1/2} = \calD(t)\quad\text{and}\quad\dom T_F^{1/2} = \calD_F.
$$
These are at the same time the domains of the closures of the forms which are induced by $T_N$, $T(t)$ and $T_F$, respectively (cf.\ \cite[Chapter VI, Theorem 2.23]{k}). Moreover (see \cite[Theorem 4.1]{hms} and \cite[Chapter VI, Theorem 2.21]{k}) for all $t\in\calI$ the following relations hold:
\begin{equation}\label{e:neumann}
\|T_N^{1/2}f\|_a\,\le\,\|T(t)^{1/2}f\|_a\quad\text{for }\,f\in\calD(t)
\end{equation}
and
\begin{equation}\label{e:friedrichs}
\|T(t)^{1/2}f\|_a\,\le\,\|T_F^{1/2}f\|_a\quad\text{for }\,f\in\calD_F.
\end{equation}
We mention that it is no restriction to assume that zero is not a turning point of $w$. Hence, the number of turning points of $w$ in $[0,a]$ is even. Let $x_1,\ldots,x_{2n}\in (0,a)$ be the turning points of $w$. Due to \cite[Section 3]{cl} there exists a uniformly positive operator $X_a\in\calB(\calH_a)$ with $X_a\calD_N\subset\calD_N$ such that for $f\in\calD_N$ we have $(X_af)(x_j) = 0$, $j=1,\ldots,2n$, and $X_af = f$ in neighborhoods of $0$ and $a$. Hence, the bounded, boundedly invertible and $J_a$-nonnegative operator $W_a := J_aX_a$ satisfies
\begin{equation}\label{e:inside1}
W_a\calD_N\subset\calD_N\quad\text{and}\quad W_a\calD(t)\subset\calD(t)\;\text{ for all }t\in\calI.
\end{equation}
By $\iota\in\{-1,1\}$ denote the sign of $w$ on $[0,x_1)\cup (x_{2n},a]$. Due to the properties of $X_a$ and \eqref{e:inside1} we have
\begin{equation}\label{e:inside2}
(W_a - \iota I)\calD_N\subset\calD_F.
\end{equation}
Now, define the operator $\wt W\in\calB(L^2(\calI,\calH_a))$ by
$$
\big(\wt WF\big)(t) := W_aF(t),\quad F\in L^2(\calI,\calH_a),\;t\in\calI.
$$
This operator is $\wt J$-nonnegative and boundedly invertible. In the following we shall show the relation
\begin{equation}\label{e:inside3}
\wt W\dom\wt T^{1/2}\,\subset\,\dom\wt T^{1/2},
\end{equation}
where $\wt T$ denotes the multiplication operator with the family $\{T(t) : t\in\calI\}$, cf.\ \eqref{e:wtT}. Clearly, $\wt T^{1/2}$ coincides with the multiplication operator with the family $\{T(t)^{1/2} : t\in\calI\}$. In particular,
$$
\dom\wt T^{1/2} = \big\{F\in L^2(\calI,\calH_a) : F(t)\in\calD(t)\;\text{a.e.},\;T(\cdot)^{1/2}F(\cdot)\in L^2(\calI,\calH_a)\big\}.
$$
In order to prove \eqref{e:inside3} let $F\in\dom\wt T^{1/2}$. Then $\wt WF\in L^2(\calI,\calH_a)$, and \eqref{e:inside1} gives $(\wt WF)(t) = W_aF(t)\in\calD(t)$ for a.e.\ $t\in\calI$. It remains to prove that $T(\cdot)^{1/2}W_aF(\cdot)$ is contained in $L^2(\calI,\calH_a)$. By the closed graph theorem and \eqref{e:inside2} there exists some $c > 0$ such that
$$
\|T_F^{1/2}(W_a - \iota I)f\|^2\,\le\,c\left( \|f\|^2 + \|T_N^{1/2}f\|^2\right),\quad f\in\calD_N.
$$
This, together with the relations \eqref{e:neumann} and \eqref{e:friedrichs}, implies
\begin{align*}
\|T(t)^{1/2}(W_a - \iota I)F(t)\|^2
&\le \|T_F^{1/2}(W_a - \iota I)F(t)\|^2\\
&\le c\left( \|F(t)\|^2 + \|T_N^{1/2}F(t)\|^2\right)\\
&\le c\Big( \|F(t)\|^2 + \|T(t)^{1/2}F(t)\|^2\Big).
\end{align*}
Hence, the function $t\mapsto T(t)^{1/2}W_aF(t) - \iota T(t)^{1/2}F(t)$ is an element of $L^2(\calI,\calH_a)$. Therefore, also $t\mapsto T(t)^{1/2}W_aF(t)$ belongs to $L^2(\calI,\calH_a)$, and \eqref{e:inside3} is proved. It is now a consequence of \eqref{e:inside3} and \cite[Proposition 3.5]{cu} that $\infty$ is not a spectral singularity of $\wt A$ and thus neither of $A$.

{\bf 2.} In the general case there exists $\eta > 0$ such that $T + \eta$ is uniformly positive. The operator $T + \eta$ is the maximal operator corresponding to
$$
\frakt_\eta(f) := \frakt(f) + \eta f = \frac 1 {|w|}\Big(\!-(pf')' + (q + \eta |w|)f\Big).
$$
Similarly, $A + \eta J = J(T + \eta)$ is the maximal operator associated with
$$
\fraka_\eta(f) := \frac 1 {w}\Big(\!-(pf')' + (q + \eta |w|)f\Big).
$$
By step 1 of this proof, $\infty$ is not a spectral singularity of $A + \eta J$. Equivalently (see  \cite[Proposition 3.5]{cu}), there exists a bounded and boundedly invertible $J$-nonnegative operator $W$ in $L^2_{|w|}(\R)$ with $W\dom A\subset\dom A$. Let $\calR > 0$ such that $[\calR,\infty)$ is of positive type, $(-\infty,-\calR]$ is of negative type with respect to $A$ and $\sigma(A)\setminus\R\subset B_\calR(0)$, and let $E$ be the spectral projection corresponding to $B_\calR(0)$, cf.\ Remark \ref{r:gsf}. Then both $E$ and $E_\perp := I - E$ are $J$-self-adjoint. Moreover, both $\ran E$ and $\ran E_\perp$ are $A$-invariant and $A_\perp := A|\ran E_\perp$ is $J$-nonnegative and boundedly invertible. Set $W_\perp := E_\perp (W|\ran E_\perp)$. Then for $f\in\ran E_\perp$ we have $[W_\perp f,f] = [Wf,f]$ which implies that also $W_\perp$ is $J$-nonnegative and boundedly invertible. In addition,
$$
W_\perp\dom A_\perp\subset E_\perp W\dom A\subset E_\perp\dom A = \dom A_\perp.
$$
Therefore, $\infty$ is not a spectral singularity of $A_\perp$ by \cite[Proposition 3.5]{cu} and thus neither of $A$.
\end{proof}

\begin{rem}
The assertion in step 1 of the proof of Theorem {\rm\ref{t:infty}} has been proved similarly but in less detail in \cite{ko}.
\end{rem}

\begin{cor}
Under the conditions on $w$ and $p$ in Theorem {\rm\ref{t:infty}} the operator $A$ is a direct sum of a bounded operator and a self-adjoint operator in a Hilbert space.
\end{cor}
\begin{proof}
Choose $\calR > 0$ as in step 2 of the proof of Theorem {\rm\ref{t:infty}} and let $E_b$ be the spectral projection of $A$ corresponding to $B_\calR(0)$. In addition, denote by $E_\pm$ the spectral projection of $A$ corresponding to $\R^\pm\setminus [-\calR,\calR]$, set $A_b := A|\ran E_b$ and $A_\pm := A|\ran E_\pm$. From $\ran E_b\subset\dom A$ it follows that $A_b$ is bounded. Since $(\calR,\infty)$ is of positive type and $(-\infty,-\calR)$ is of negative type with respect to $A$, the inner product spaces $(\ran E_+,\product)$ and $(\ran E_-,-\product)$ are Hilbert spaces. Moreover, $\ran(I - E_b) = \ran E_+[\ds]\ran E_-$. Therefore, $A_s := A_+[\ds]A_-$ is self-adjoint in the Hilbert space $(\ran(I - E_b),\langle\cdot,\cdot\rangle)$, where $\langle f_+ + f_-,g_+ + g_-\rangle = [f_+,g_+] - [f_-,g_-]$, $f_\pm,g_\pm\in\ran E_\pm$, and $A = A_b[\ds]A_s$.
\end{proof}

\section*{Acknowledgements}
The author thanks Aleksey Kostenko and Fritz Gesztesy for fruitful discussions in Voronezh and Columbia, MO, respectively. Moreover, he is grateful to Jussi Behrndt and Carsten Trunk for their constant support.


\section*{Contact information}
Friedrich Philipp: Technische Universit\"at Berlin, Institut f\"ur Mathematik, Stra\ss e des 17.\ Juni 136, 10623 Berlin, Germany, philipp@math.tu-berlin.de
\end{document}